\documentclass[10pt]{article}
\usepackage[utf8]{inputenc}
\usepackage{bussproofs}
\usepackage{mathrsfs}
\usepackage{amsmath}
\usepackage{amsthm}
\usepackage{amssymb}
\usepackage[dvipsnames]{xcolor}
\usepackage{stmaryrd}
\usepackage{calc}
\usepackage{cmll}
\usepackage{url}
\usepackage{amsmath}
\usepackage{longtable}
\usepackage{stmaryrd}
\usepackage{{xargs,ifthen,xstring}}
\usepackage{authblk}
\usepackage{txfonts}
\usepackage{array}
\newtheorem{theorem}{Theorem}

\newtheorem{definition}{Definition}
\newtheorem{proposition}{Proposition}

\newtheorem{corollary}{Corollary}
\newtheorem*{fact}{Fact}

\title{Some results in non-monotonic\\proof-theoretic semantics}

\author{Antonio Piccolomini d'Aragona\\

University of T\"{u}bingen, T\"{u}bingen, Germany\\

\texttt{antonio.piccolomini-daragona@uni-tuebingen.de}}
\date{}

\begin{document}

\maketitle

\begin{abstract}
   I explore the relationships between Prawitz's approach to non-monotonic proof-theoretic validity, which I call \emph{reducibility semantics}, and a later proof-theoretic approach, that I call \emph{standard base semantics}. I prove that, if suitable conditions are met, reducibility semantics and standard base semantics are equivalent. As a side-result, I show that a similar relation holds, albeit in a weaker way, also between reducibility semantics and a variant of standard base semantics due to Sandqvist. Finally, notions of ‘‘point-wise" soundness and completeness (called base-soundness and base-completeness) are discussed against certain known principles from the proof-theoretic literature, as well as against monotonic proof-theoretic semantics. Intuitionistic logic is proved not to be ‘‘point-wise" complete on any kind of non-monotonic proof-theoretic semantics. The way in which this result is proved, as well as the overall behaviour of ‘‘point-wise" soundness and completeness, is significantly different in the non-monotonic framework as compared to what happens in the monotonic one---where the notions at issue can be used too to prove the ‘‘point-wise" incompleteness of intuitionistic logic.
\end{abstract}

\paragraph{Keywords} Proof-theoretic semantics, non-monotonicity, consequence, completeness

\section{Introduction}

Proof-theoretic semantics (henceforth PTS) is a formal, proof-based semantics stemming from Prawitz's normalisation results for Gentzen's Natural Deduction \cite{francez, gentzenuntersuchungen, prawitz1965, schroeder-heisterSE}. In its original formulation, due to Prawitz himself \cite{prawitz1971, prawitz1973}, PTS is based on argument structures and reductions for non-introduction inferences, which are in turn used for defining valid arguments, i.e., linguistic representations of proofs. Argumental validity is relativised to theories governing deduction at the atomic level, and it is later on generalised by universally quantifying over all such theories, so as to account for logical argumental validity. The notion of (logical) consequence is a derived one, as it is understood in terms of existence of a (logically) valid argument from suitable assumptions to a suitable conclusion. I will call this kind of PTS \emph{reducibility semantics}.

In the current mainstream approach, (logical) consequence is defined outright instead, so that argument structures and reductions are left out. The constructivist burden is all put on atomic theories, as the latter are understood (also by Prawitz) as fixing the meaning of the non-logical terminology in terms of the inferential role played by this terminology in atomic deduction. This type of PTS, which I shall refer to as \emph{standard base semantics}, has been investigated extensively in recent years, while a variant of it due to Sandqvist \cite{sandqvist}---whence \emph{Sandqvist's base semantics}---deals with disjunction in an ‘‘elimination-like" fashion, i.e., it explains the validity of a disjunction not in terms of validity of either disjunct, but in terms of categorical validity of atoms under hypothetical validity of atoms under both disjuncts.

All the three approaches to PTS mentioned so far can be developed mainly in two ways: \emph{monotonically} or \emph{non-monotonically}. The distinction boils down to how one deals with the ‘‘open" case, that is, validity of an argument with non-empty set of assumptions (in the case of reducibility semantics), or consequence of $A$ from $\Gamma \neq \emptyset$ (in the case of the two base semantics). Monotonic approaches adopt a closure condition throughout arbitrary extensions of an underlying atomic theory while, in the non-monotonic framework, the reference to extensions is left out, and the closure condition always operates on one and the same atomic theory.

Many completeness and incompleteness results have been proved for the monotonic version of both base semantics and its Sandqvist's variant \cite{piechadecampossanz, piechaschroeder-heisterdecampossanz, piechaschroeder-heister2019, sandqvist, schroederheisterrolf, stafford1, stafford2}. In \cite{piccolominiinversion}, I investigated some relationships between the different kinds of monotonic PTS. I also proved some equivalence and connection results, and applied them to issues of completeness of logics over these semantics. This led me to define notions of ‘‘point-wise" soundness and completeness, called base-soundness and base-completeness, and to discuss some of their properties against the kinds of PTS at issue here, and against the results one could prove about them, or about their relation. This paper is a follow-up where I aim at showing how and to what extent the results presented in \cite{piccolominiinversion} can be applied to the non-monotonic context. I think this is a relevant topic, for at least two reasons.

First, Prawitz has used a monotonic approach only in \cite{prawitz1971}, and the non-monotonic version instead from then on. Hence, he seems to prefer non-monotonic PTS over the monotonic variant. The ultimate reason for this is that, as said, the atomic theories which argumental validity and consequence are defined over, are also taken as ‘‘models", i.e., as systems which fix the meaning of the non-logical terminology by showing how this terminology behaves at the level of atomic derivability. Contrarily to the Tarskian idea that non-logical meaning is determined by mappings from the language onto (mostly set-theoretic) structures, this is more consonant with broad constructivist desiderata. For example, it can be seen as nothing but a full-application of Gentzen's well-known idea that the meaning of the logical constants is defined by special rules for these constants (specifically, introduction rules) \cite{gentzen}: that the meaning of the non-logical dictionary stems from the deductive behaviour of this dictionary at the atomic level, just means that this meaning is given by atomic rules which pertain to those non-logical signs. More generally, the idea is also in line with Dummett's requirement for a verificationist theory where meaning determination is not detached from linguistic practice \cite{dummett1991}, hence with (an interpretation of) Wittgenstein's slogan that ‘‘meaning is use" \cite{wittgenstein}. This said, it is easy to see what might seem problematic in Prawitz's view when a monotonic approach is endorsed: to require that an argument should be valid or that a consequence should hold over arbitrary extensions of a given atomic theory, may not completely capture the idea that that validity and consequence depend on a \emph{fixed} meaning of the non-logical terminology. For, as soon as we change the base, we also change the meaning of this terminology.

The second reason why non-monotonic PTS is of interest has to do with how advocates of the \emph{monotonic} approach normally defend their position. They remark that it is part of our understanding of argumental validity and of deductive consequence that these notions are stable over extensions of one's knowledge base---on this, see e.g. \cite{piechaschroeder-heisterbases}. We can certainly agree on that, while also remarking, however, that this requires changing the way we look at (the hierarchy of) atomic theories: they are no longer (or not only) ‘‘models", but inference-networks representing one's epistemic states. From this point of view, monotonic and non-monotonic PTS are nothing but sides of the same coin, in the sense that they capture complementary intuitions about argumental validity and deductive consequence. Roughly put, a proof can be said to be epistemically convincing because of \emph{two} components at least: the meaning of the sentential and sub-sentential expressions it involves, and the validity of the inferences of which it is made. Depending on whether one takes one or the other aspect to be more prominent, one will have a more semantically or, alternatively, a more deductively oriented perspective on proof-based approaches to meaning. Monotonic PTS is more deductively oriented, whereas non-monotonic PTS constitutes a semantic counter-part where it is very natural to say that an argument is valid or that a consequence holds relative to a given atomic theory, while failing on extensions of this theory---just like, for example, the induction principle holds only because of the special structure of the set of the natural numbers (i.e., because of the rules governing the numerical signs), and is thus not expected to hold on all the extensions of this set, e.g., on the real numbers.\footnote{This situation is not random, as it can be shown to stem from some formal and conceptual symmetries that induce an order relation among potential approaches to PTS---I have dealt with this in \cite{ptssquare}.}

There is however an independent interest---which abstracts from the monotonic/non-monotonic split---in comparing reducibility semantics, on the one hand, with standard base semantics and Sandqvist's semantics, on the other. This articulates into two points again.

The first point is mainly about the relationship between reducibility semantics and standard base semantics. As said, in reducibility semantics the notion of (logical) consequence is defined in a derivative way, based on a prior notion of argumental validity. Thus, the (logical) consequence relation is always ‘‘witnessed" by (the existence of) a suitable argument structure, equipped with suitable reductions for non-introduction inferences. In standard base semantics, instead, argument structures and reductions are left out, so it becomes very natural to wonder whether the removal of the (logical) consequence ‘‘witness" also implies some loss in constructivity. This is not to suspect, of course, that standard base semantics \emph{is not} constructive, as constructivity is in a sense secured by the special kind of ‘‘models" which this approach shares with reducibility semantics. However, the fact that the holding of a (logical) consequence relation is not grounded on a valid argument ‘‘witnessing" that the relation does hold, might make one suspicious about the fact that standard base semantics is in a way ‘‘less" constructive than reducibility semantics. I shall show below that, to a large extent, and provided that certain conditions are met, this is in fact not the case.

What I just said applies to the comparison of reducibility semantics and Sandqvist's base semantics as well, for the latter does without argument structures and reductions too. But a second point of interest enters the scene here. In \emph{monotonic} PTS, intuitionistic logic can be proved to be complete over Sandqvist's approach (with ‘‘models" of a special kind), and incomplete over standard base semantics and reducibility semantics (independently of the kind of ‘‘models", but including those where Sandqvist's completeness obtains). Now, given that the ‘‘models" underlying both approaches overlap, the completeness/incompleteness mismatch makes it worth investigating whether the two frameworks are comparable at the level of consequence over a ‘‘model".

The main aim of this paper is that of comparing non-monotonic reducibility semantics with standard non-monotonic base semantics. The main result is, as said, that the former is equivalent to the latter in a significant---albeit not most general---sense. On the other hand, since no (in)completeness theorem is currently available for non-monotonic Sandqvist's base semantics, the interest of comparing the latter with non-monotonic reducibility semantics is in a way less relevant. For the reasons mentioned above, however, I shall very quickly mention some results, to the effect that these two approaches are comparable only in a ‘‘conditionally global" way.

A second aim of the paper is comparing monotonic and non-monotonic PTS relative to the comparability of the kinds of PTS at issue here. All the comparability results proved for the monotonic picture hold, \emph{mutatis mutandis}, in the non-monotonic picture as well. It is therefore natural to wonder whether this applies also to the \emph{consequences} of these results. The interaction of the latter with the aforementioned notions of base-soundness and base-completeness, has some quite interesting implications in the monotonic context, providing for example sufficient conditions for reducibility semantics and standard base semantics to be equivalent to Sandqvist's semantics, as well as for given logics to be complete over reducibility semantics and standard base semantics. These conditions are however of little use with intuitionistic logic (and a number of other logics too), essentially because the latter can be proved to be base-incomplete (with respect to ‘‘models" of any kind) over any of the PTS approaches under consideration here. Base-incompleteness of intuitionistic logic also holds in the non-monotonic approach but, due to non-monotonicity, the reasons for this lie in broader results of a different kind than those used for proving the same in monotonic PTS. Thus, the discussion on base-soundness and base-completeness leads to a further appreciation of the differences between monotonic and non-monotonic PTS.

The paper is structured as follows. In Section 2, I introduce atomic theories which reducibility semantics and base semantics are to be defined over. In Section 3, I define the three semantics themselves (in a non-monotonic way), and prove basic facts about them. In Section 4, I prove results which provide a comparison of the three approaches (in the sense hinted at above). In Section 5, I introduce notions of ‘‘point-wise" soundness and completeness, and investigate their behaviour as compared to that of their monotonic counter-parts. As said, in Sections 3 and 4 the focus on Sandqvist's variant will be very quick, and only due to the concern of having an as much as possible complete picture of how the major existing approaches to non-monotonic proof-theoretic semantics interact with each other.

\section{Language and atomic bases}

I will restrict myself to a propositional language, defined as follows.

\begin{definition}
The \emph{language} $\mathscr{L}$ is given by the grammar:

\begin{center}
    $X := p, q, r, s, t, ... \ | \ \bot \ | \ X \wedge X \ | \ X \vee X \ | \ X \rightarrow X$
\end{center}
where $\bot$ is an atomic constant symbol for absurdity and $\neg X := X \rightarrow \bot$.
\end{definition}
\noindent I shall indicate the set of the atomic formulas and the set of the formulas of $\mathscr{L}$ respectively by $\texttt{ATOM}_{\mathscr{L}}$ and $\texttt{FORM}_{\mathscr{L}}$. I shall use capital Latin letters $A, B, C, ...$ for formulas, and capital Greek letters $\Gamma, \Delta, \Theta, ...$ for sets of formulas.

\begin{definition}
Any $A \in \texttt{ATOM}_{\mathscr{L}}$ is an \emph{atomic rule of level} $0$. Given atomic rules whose maximal level is $k - 1$
    \begin{prooftree}
        \AxiomC{$\mathfrak{C}_1$}
        \noLine
        \UnaryInfC{$A_1$}
        \AxiomC{}
        \noLine
        \UnaryInfC{$\dots$}
        \AxiomC{$\mathfrak{C}_n$}
        \noLine
        \UnaryInfC{$A_n$}
        \noLine
        \TrinaryInfC{}
    \end{prooftree}
where each $\mathfrak{C}_i$ may be empty \emph{(}$i \leq n$\emph{)},
    \begin{prooftree}
        \AxiomC{$[\mathfrak{C}_1]$}
        \noLine
        \UnaryInfC{$A_1$}
        \AxiomC{}
        \noLine
        \UnaryInfC{$\dots$}
        \AxiomC{$[\mathfrak{C}_n]$}
        \noLine
        \UnaryInfC{$A_n$}
        \TrinaryInfC{$A$}
    \end{prooftree}
with $A \in \texttt{ATOM}_{\mathscr{L}}$ is an atomic rule of level $k$.
\end{definition}
\noindent As usual, square brackets indicate that the rule discharges the corresponding assumptions. However, what is assumed (and thus discharged) here are not atomic formulas only, but also atomic rules of a lower level than the rule which discharges them---i.e., an atomic rule of level $k$ discharges atomic rules of level $\leq k - 2$, see \cite{piechaschroeder-heisterdecampossanz, piechaschroeder-heisterbases, piechaschroeder-heister2019}.

\begin{definition}
    An \emph{atomic base of level} $k$ is a set of atomic rules of maximal level $k$.
\end{definition}
\noindent Atomic derivations are defined by induction in a usual way. When undischarged rules of an atomic derivation with conclusion $A$ belong to a given base or to a set of atomic rules $\mathfrak{C}$, $A$ is said to be derivable from $\mathfrak{C}$ in the base---for more precise definitions, see \cite{piechaschroeder-heisterdecampossanz}. The derivability relation over $\mathfrak{B}$ is indicated by $\vdash_\mathfrak{B}$, while the derivation set of $\mathfrak{B}$ is indicated by $\texttt{DER}_\mathfrak{B}$. Atomic bases are assumed to be always consistent, i.e., $\forall \mathfrak{B} \ (\not\vdash_\mathfrak{B} \bot)$.

\section{Non-monotonic proof-theoretic semantics}

I shall not follow a chronological order, but start from base semantics instead. This is because, as we shall see, the equivalence and connection results that I will establish between reducibility semantics and base semantics are mostly oriented towards proving how (in)completeness results proved \emph{for the latter} can be applied \emph{to the former}.

In what follows, $\Gamma$ is always understood as a \emph{finite} set. This is to avoid some compactness issues that emerge when comparing reducibility semantics and base semantics---pointed out by \cite{stafford1}---and is sufficient for raising my points and prove a number of results.

\subsection{Base-semantics}

Consequence over an atomic base in standard base semantics is defined as follows.

\begin{definition}
    $A$ is a \emph{consequence of} $\Gamma$ \emph{on} $\mathfrak{B}$, written $\Gamma \models_\mathfrak{B} A \Longleftrightarrow$
    \begin{itemize}
    \item $\Gamma = \emptyset \Longrightarrow$
    \begin{itemize}
        \item $A \in \texttt{ATOM}_\mathscr{L} \Longrightarrow \ \vdash_\mathfrak{B} A$
        \item $A = B \wedge C \Longrightarrow \ \models_\mathfrak{B} B$ and $\models_\mathfrak{B} C$
        \item $A = B \vee C \Longrightarrow \ \models_\mathfrak{B} B$ or $\models_\mathfrak{B} C$
        \item $A = B \rightarrow C \Longrightarrow B \models_\mathfrak{B} C$
    \end{itemize}
    \item $\Gamma \neq \emptyset \Longrightarrow \ (\models_\mathfrak{B} \Gamma \Longrightarrow \ \models_\mathfrak{B} A)$
    \end{itemize}
    where $\models_\mathfrak{B} \Gamma$ means $\forall B \in \Gamma \ (\models_\mathfrak{B} B)$.
\end{definition}
\noindent Sandqvist's (non-monotonic) version, written $\models^s_\mathfrak{B}$, is defined in the same way as $\models_\mathfrak{B}$, except for the clause for disjunction.

\begin{definition}
    $\models^s_\mathfrak{B} A \vee B \Longleftrightarrow \forall C \in \texttt{ATOM}_\mathscr{L} \ (A \models^s_\mathfrak{B} C$ and $B \models^s_\mathfrak{B} C \Longrightarrow \ \models^s_\mathfrak{B} D)$.
\end{definition}
\noindent Let $\Vdash$ be either $\models$ or $\models^s$. $\Vdash$ is non-monotonic in the sense that there are $\Gamma, A$ and $\mathfrak{B}$ such that $\Gamma \Vdash_\mathfrak{B} A$ but, for some extension $\mathfrak{C}$ of $\mathfrak{B}$, $\Gamma \not\Vdash_\mathfrak{C} A$---where an extension of an atomic base is a super-set of the base. E.g., take $\Gamma = \{p_1\}, A = p_2, \mathfrak{B} = \emptyset$ and $\mathfrak{C} = \{p_1\}$.

\begin{definition}
    $\Gamma \Vdash A \Longleftrightarrow \forall \mathfrak{B} \ (\Gamma \Vdash_\mathfrak{B} A)$.
\end{definition}

\noindent Intuitionistic logic (henceforth $\texttt{IL}$) can be shown to be incomplete over $\models$ with classical meta-logic---in fact classical logic is sound and complete over $\models$ under this condition \cite{piccolomininote, piechaschroeder-heister2019, schroederheisterrolf}.\footnote{That the meta-logic needs to be classical does not imply that the result is less problematic for the project of building a semantics for $\texttt{IL}$. As Piecha and Schroeder-Heister observe, ‘‘as a negative result, [this] is as devastating for [Prawitz's conjecture of completeness of $\texttt{IL}$] as a constructive proof" \cite[244]{piechaschroeder-heister2019}. Or again, ‘‘claiming that completeness can nonetheless be proved by intuitionistic (and thus also by classical) means implies claiming a classical contradiction. Given that these proofs can be coded in first-order arithmetic and that classical arithmetic and Heyting arithmetic are equiconsistent, such a claim cannot be upheld" \cite[501]{schroederheisterrolf}.} Sandqvist's variant is instead interesting because, in the monotonic framework, $\texttt{IL}$ is complete on it for atomic bases of level $\geq 2$ \cite{sandqvist}---more on this in Section 4.3. Sandqvist's picture may play a role in the non-monotonic format too since, as shown in \cite{schroederheisterrolf}, it is an approach which prioritises elimination rules---on elimination-based approaches see also \cite{hermogenespragmatist, takemura}.\footnote{It is worth noting that, in Sandqvist's paper \cite{sandqvist}, $\bot$ is not an atomic constant symbol for absurdity, but a $0$-ary propositional connective explained by a clause which states that $\bot$ holds on an atomic base iff every atom holds on the base. This approach can be equated with the one where $\bot$ is an atomic constant symbol by postulating that atomic bases always contain an ‘‘atomic explosion rule". In \emph{monotonic} PTS, this issue plays a crucial role, as there intuitionistic logic is complete over Sandqvist's variant (with atomic bases of a special kind), and incomplete over the standard one. When comparing the different kinds of PTS relative to (in)completeness issues, we must hence find a common ground for the comparison to go through without problems---see \cite{piccolominiinversion}. For what is at stake in this paper, though, this aspect seems to me to be less important. In what follows I shall therefore abstract from this issue. I shall come back to this in Section 4.3.} But since no (in)completeness result are currently known for the non-monotonic variant of Sandqvist's base semantics, in what follows I shall deal only quickly with this variant.

\subsection{Reducibility semantics}

As said, in reducibility semantics we do not start from outright notions of consequence between $\Gamma$ and $A$, but from the notion of (logically) valid argument. Thus, some preliminary notions are required, in particular, those of argument structure and reduction for non-introduction rule. The need of defining argumental validity first, and based on that the notion of consequence, as well as the need of introducing a number of preliminary notions, explain why this section is going to be longer than that of base semantics. 

\begin{definition}
    An \emph{argument structure} is a pair $\langle T, f \rangle$ where:
    
    \begin{itemize} 
    \item $T$ is a finite rooted tree with order relation $\omega$, nodes labelled by formulas of $\mathscr{L}$, and top-nodes partitioned into two groups, i.e., axiomatic and non-axiomatic;
    \item $f$ is a function mapping onto lower nodes elements of: \emph{(}a\emph{)} a sub-set of the non-axiomatic top-nodes of $T$; \emph{(}b\emph{)} a sub-set of the axiomatic top-nodes of $T$ labelled by atoms; \emph{(}c\emph{)} a sub-set of $\wp(\omega)$, all the elements of which contain only the pairs linking a node labelled by an atom to all its children, labelled by atoms too, with no non-axiomatic top-node of $T$ mapped by $f$ onto such a maximal node, or onto the image of the latter through $f$.
    \end{itemize}
    The non-axiomatic top-nodes and the root of $T$ are called, respectively, the \emph{assumptions} and the \emph{conclusion} of $\langle T, f \rangle$.
\end{definition}
\noindent It should be clear that $f$ indicates where certain assumptions, atomic axioms, or atomic rules are discharged in $T$.

\begin{definition}
    Given $\mathscr{D} = \langle T, f \rangle$, we say that $\mathscr{D}$ is \emph{closed} if the domain of $f$ contains all the non-axiomatic top-nodes of $T$, otherwise we say that $\mathscr{D}$ is \emph{open}. When $\mathscr{D}$ has undischarged assumptions $\Gamma$ and conclusion $A$, we say it is an argument structure \emph{from} $\Gamma$ \emph{to} or \emph{for} $A$.
\end{definition}
\noindent As usual, an argument structure $\mathscr{D}$ from $\Gamma$ to $A$ is also indicated by the figure

\begin{prooftree}
    \AxiomC{$\Gamma$}
    \noLine
    \UnaryInfC{$\mathscr{D}$}
    \noLine
    \UnaryInfC{$A$}
\end{prooftree}
The notion of (\emph{immediate}) \emph{sub-structure} of an argument structure can be easily defined. The same applies to the notion of \emph{replacement} of a sub-structure $\mathscr{D}^*$ of an argument structure $\mathscr{D}$ with an argument structure $\mathscr{D}^{**}$---written $\mathscr{D}[\mathscr{D}^{**}/\mathscr{D}^*]$---although in this case one may have to take care of potential re-indexing of the discharge functions of $\mathscr{D}$ and $\mathscr{D}^{**}$. I will abstract from these details here.

\begin{definition}
    Given $\mathscr{D}$ from $\Gamma = \{A_1, ..., A_n\}$ to $A$, and a function $\sigma$ from formulas to argument structures such that $\sigma(A_i)$ is a \emph{(}closed\emph{)} argument structure for $A_i$ $(i \leq n)$, we say that $\mathscr{D}^\sigma = \mathscr{D}[\sigma(A_1), ..., \sigma(A_n)/A_1, ..., A_n]$ is a \emph{(closed)} \emph{instance} of $\mathscr{D}$.
\end{definition}
\noindent Let me now give an example of argument structure. Axiomatic top-nodes will be distinguished from non-axiomatic top-nodes by putting an horizontal bar on top of the former.

\begin{prooftree}
    \AxiomC{$[p \wedge \neg q]_1$}
    \AxiomC{}
    \UnaryInfC{$\neg r$}
    \BinaryInfC{$q$}
    \AxiomC{$s$}
    \RightLabel{$[R]_2$}
    \BinaryInfC{$t$}
    \AxiomC{$\overline{[q]}_3$}
    \AxiomC{$\neg \neg p \vee r$}
    \BinaryInfC{$q \rightarrow \neg s$}
    \BinaryInfC{$r$}
    \RightLabel{$3$}
    \UnaryInfC{$t$}
    \AxiomC{}
    \UnaryInfC{$s$}
    \RightLabel{$2$}
    \BinaryInfC{$p$}
    \UnaryInfC{$q \rightarrow (p \vee r)$}
    \RightLabel{$1$}
    \UnaryInfC{$p \wedge \neg q \rightarrow (q \rightarrow (p \vee r))$}
\end{prooftree}
This is an open argument structure, from assumptions $s$ and $\neg \neg p \vee r$ to $p \wedge \neg q \rightarrow (q \rightarrow (p \vee r))$. The discharge function maps the non-axiomatic top-node labelled by $p \wedge \neg q$ onto the root-node, the axiomatic top-nodes labelled by $q$ onto the node labelled by $t$ with upper edge labelled by $3$, and the set of edges labelled by $R$ onto the node labelled by $p$ with upper edges labelled by $2$. The structure is an open instance of itself, while a closed instance is the following.

\begin{prooftree}
    \AxiomC{$[p \wedge \neg q]_1$}
    \AxiomC{}
    \UnaryInfC{$\neg r$}
    \BinaryInfC{$q$}
    \AxiomC{}
    \UnaryInfC{$p \vee q$}
    \UnaryInfC{$s$}
    \RightLabel{$[R]_2$}
    \BinaryInfC{$t$}
    \AxiomC{$\overline{[q]}_3$}
    \AxiomC{}
    \UnaryInfC{$p$}
    \AxiomC{$[q \rightarrow s]_4$}
    \BinaryInfC{$q \vee p$}
    \RightLabel{$4$}
    \UnaryInfC{$\neg \neg p \vee r$}
    \BinaryInfC{$q \rightarrow \neg s$}
    \BinaryInfC{$r$}
    \RightLabel{$3$}
    \UnaryInfC{$t$}
    \AxiomC{}
    \UnaryInfC{$s$}
    \RightLabel{$2$}
    \BinaryInfC{$p$}
    \UnaryInfC{$q \rightarrow (p \vee r)$}
    \RightLabel{$1$}
    \UnaryInfC{$p \wedge \neg q \rightarrow (q \rightarrow (p \vee r))$}
\end{prooftree}

\begin{definition}
    An \emph{inference} is a triple $\langle \langle \mathscr{D}_1, ..., \mathscr{D}_n \rangle, A, \delta \rangle$ $(n \geq 0)$ where $\delta$ is an extension of the discharge functions associated to the $\mathscr{D}_i$-s $(i \leq n)$. An \emph{inference rule} is a set of inferences.
\end{definition}
\noindent Based on the previous definition, it is easy to see that an inference $\langle \langle \mathscr{D}_1, ..., \mathscr{D}_n \rangle, A, \delta \rangle$ is uniquely associated to an argument structure
    \begin{prooftree}
        \AxiomC{$\mathscr{D}_1, ..., \mathscr{D}_n$}
        \RightLabel{$\delta$}
        \UnaryInfC{$A$}
    \end{prooftree}
obtained by conjoining the trees of the $\mathscr{D}_i$-s via a root-node $A$, and by adding $\delta$ to the discharge functions of the $\mathscr{D}_i$-s $(i \leq n)$. As regards inference rules, I shall assume that they can be always identified via some meta-linguistic descriptions. Thus, for example, the standard $\wedge$-elimination and $\vee$-elimination in Gentzen's natural deduction, are

\begin{center}
    $\wedge_{E, i} = \{\langle \mathscr{D}, A_i, \emptyset \rangle \ | \ \mathscr{D}$ for $A_1 \wedge A_2, i = 1, 2\}$
\end{center}

\begin{center}
    $\vee_E = \{\langle \langle \mathscr{D}_1, \mathscr{D}_2, \mathscr{D}_3 \rangle, C, \{f(A_{\mathscr{D}_2}) = C$, $f(B_{\mathscr{D}_3}) = C\} \rangle \ |$
    
    $\mathscr{D}_1$ for $A \vee B, \mathscr{D}_2$ from $A_n$ to $C, \mathscr{D}_3$ from $B_m$ to $C\}$
\end{center}
where $A_{\mathscr{D}_2}$ and $B_{\mathscr{D}_3}$ are the occurrences of $A$ and $B$ in $\mathscr{D}_2$ and $\mathscr{D}_3$ respectively. We may consider two additional rules.

\begin{center}
    $\texttt{Wk} = \{\langle \mathscr{D}, (A \wedge C) \rightarrow B, \emptyset \rangle \ | \ \mathscr{D}$ for $A \rightarrow B\}$
\end{center}

\begin{center}
    $\vee_\lambda = \{\langle \mathscr{D}, A, \emptyset \rangle \ | \ \mathscr{D}$ for $A \vee B\}$
\end{center}
\noindent All the rules above can be described by the following inference schemata, respectively.

\begin{prooftree}
    \AxiomC{$A_1 \wedge A_2$}
    \RightLabel{$\wedge_{E, i}$}
    \UnaryInfC{$A_i$}
    \AxiomC{$A \vee B$}
    \AxiomC{$[A]$}
    \noLine
    \UnaryInfC{$C$}
    \AxiomC{$[B]$}
    \noLine
    \UnaryInfC{$C$}
    \RightLabel{$\vee_E$}
    \TrinaryInfC{$C$}
    \AxiomC{$A \rightarrow B$}
    \RightLabel{$\texttt{Wk}$}
    \UnaryInfC{$(A \wedge C) \rightarrow B$}
    \AxiomC{$A \vee B$}
    \RightLabel{$\vee_\lambda$}
    \UnaryInfC{$A$}
    \noLine
    \QuaternaryInfC{}
\end{prooftree}

\begin{definition}
    $\mathscr{D}$ is \emph{canonical} iff it is associated to an instance of an introduction rule, while it is \emph{non-canonical} otherwise.
\end{definition}
\noindent So, the argument structure exemplified above is canonical, with immediate sub-structure non-canonical.

\begin{definition}
    Given a rule $R$, a \emph{reduction} for $R$ is a function $\phi$ from and to argument structures, defined on some sub-set $\mathbb{D}$ of $R$ such that, for every $\mathscr{D}$ associated to any element of $\mathbb{D}$ and every $\sigma$
    \begin{itemize}
        \item[1.] the instance of $R$ which $\mathscr{D}^\sigma$ is associated to, is also in $\mathbb{D}$
        \item[2.] $\phi$ is defined on $\mathscr{D}^\sigma$, and $\phi(\mathscr{D}^\sigma) = \phi(\mathscr{D})^\sigma$
        \item[3.] $\mathscr{D}$ is from $\Gamma$ to $A \Longrightarrow \phi(\mathscr{D})$ is from $\Gamma^* \subseteq \Gamma$ to $A$.
    \end{itemize}
\end{definition}
\noindent Hence, for example, the standard reductions $\phi_\wedge, \phi_\vee$ for $\wedge$- and $\vee$-detours are defined on the following respective sets.

\begin{center}
    $\{\langle \mathscr{D}, A_i, \emptyset \rangle \ | \ \mathscr{D}$ canonical$\} \subset \wedge_{E, i}$
\end{center}

\begin{center}
    $\{\langle \langle \mathscr{D}_1, \mathscr{D}_2, \mathscr{D}_3 \rangle, C, \delta \rangle \ | \ \mathscr{D}_1$ canonical$\} \subset \vee_E$
\end{center}
They behave in the usual way.

\begin{prooftree}
    \AxiomC{$\mathscr{D}_1$}
    \noLine
    \UnaryInfC{$A_1$}
    \AxiomC{$\mathscr{D}_2$}
    \noLine
    \UnaryInfC{$A_2$}
    \BinaryInfC{$A_1 \wedge A_2$}
    \RightLabel{$\wedge_{E, i}$}
    \UnaryInfC{$A_i$}
    \AxiomC{$\stackrel{\phi_\wedge}{\Longrightarrow}$}
    \noLine
    \UnaryInfC{}
    \AxiomC{$\mathscr{D}_i$}
    \noLine
    \UnaryInfC{$A_i$}
    \noLine
    \TrinaryInfC{}
\end{prooftree}

\begin{prooftree}
    \AxiomC{$\mathscr{D}_1$}
    \noLine
    \UnaryInfC{$A_i$}
    \UnaryInfC{$A_1 \vee A_2$}
    \AxiomC{$[A_1]_1$}
    \noLine
    \UnaryInfC{$\mathscr{D}_2$}
    \noLine
    \UnaryInfC{$C$}
    \AxiomC{$[A_2]_2$}
    \noLine
    \UnaryInfC{$\mathscr{D}_3$}
    \noLine
    \UnaryInfC{$C$}
    \RightLabel{$\vee_E, 1, 2$}
    \TrinaryInfC{$C$}
    \AxiomC{$\stackrel{\phi_\vee}{\Longrightarrow}$}
    \noLine
    \UnaryInfC{}
    \AxiomC{$\mathscr{D}_1$}
    \noLine
    \UnaryInfC{$[A_i]$}
    \noLine
    \UnaryInfC{$\mathscr{D}_{i + 1}$}
    \noLine
    \UnaryInfC{$C$}
    \noLine
    \TrinaryInfC{}
\end{prooftree}
\noindent We may also consider reductions $\phi_{\texttt{Wk}}$ and $\phi_{\vee_\lambda}$ for $\texttt{Wk}$ and $\vee_\lambda$ defined on the following respective sets.

\begin{center}
    $\{\langle \mathscr{D}, (A \wedge C) \rightarrow B, \emptyset \rangle \ | \ \mathscr{D}$ canonical$\} \subset \texttt{Wk}$
\end{center}

\begin{center}
    $\{\langle \mathscr{D}, A, \emptyset \rangle \ | \ \mathscr{D}$ canonical with immediate sub-structure for $A\} \subset \vee_\lambda$
\end{center}
The reductions behave in the following respective ways.

\begin{prooftree}
    \AxiomC{$[A]_1$}
    \noLine
    \UnaryInfC{$\mathscr{D}^*$}
    \noLine
    \UnaryInfC{$B$}
    \RightLabel{$1$}
    \UnaryInfC{$A \rightarrow B$}
    \RightLabel{$\texttt{Wk}$}
    \UnaryInfC{$(A \wedge C) \rightarrow B$}
    \AxiomC{$\stackrel{\phi_{\texttt{Wk}}}{\Longrightarrow}$}
    \noLine
    \UnaryInfC{}
    \AxiomC{$[A \wedge C]_1$}
    \RightLabel{$\wedge_{E, 1}$}
    \UnaryInfC{$A$}
    \noLine
    \UnaryInfC{$\mathscr{D}^*$}
    \noLine
    \UnaryInfC{$B$}
    \RightLabel{$1$}
    \UnaryInfC{$(A \wedge C) \rightarrow B$}
    \noLine
    \TrinaryInfC{}
\end{prooftree}

\begin{prooftree}
    \AxiomC{$\mathscr{D}^*$}
    \noLine
    \UnaryInfC{$A$}
    \UnaryInfC{$A \vee B$}
    \RightLabel{$\vee_\lambda$}
    \UnaryInfC{$A$}
    \AxiomC{$\stackrel{\phi_{\vee_\lambda}}{\Longrightarrow}$}
    \noLine
    \UnaryInfC{}
    \AxiomC{$\mathscr{D}^*$}
    \noLine
    \UnaryInfC{$A$}
    \noLine
    \TrinaryInfC{}
\end{prooftree}

\noindent A set of reductions will be indicated by $\mathfrak{J}$. An \emph{extension} of $\mathfrak{J}$ is any set of reductions $\mathfrak{J}^+$ such that $\mathfrak{J} \subseteq \mathfrak{J}^+$. In what follows, I shall not distinguish simple extensions, from what \cite{prawitz1973} calls \emph{consistent} extensions of sets of reductions. In other words, we shall abstract from whether we are allowing or not for what \cite{schroeder-heister2006} calls \emph{alternative justifications} for argument structures.\footnote{A consistent extension of $\mathfrak{J}$ is any $\mathfrak{J}^+ \supseteq \mathfrak{J}$ such that, for every $\phi \in \mathfrak{J}^+ - \mathfrak{J}$, if $\phi$ is a reduction for $R$, then no $\phi^* \in \mathfrak{J}$ is a reduction for $R$.} Nor shall I put restrictions on the kind of functions which (sets of) reductions may be, in particular, I shall not \emph{a priori} require them to be necessarily \emph{constructive} functions of a schematic kind, as implicitly done in \cite{prawitz1973}---but observe that the reductions defined above \emph{are} constructive. I will quickly come back to this point in the concluding remarks---for a preliminary address, see \cite{piccolomininote}.

\begin{definition}
    $\mathscr{D}$ \emph{immediately reduces} to $\mathscr{D}^*$ \emph{via} $\mathfrak{J}$ iff, for some sub-structure $\mathscr{D}^{**}$ of $\mathscr{D}$, there is $\phi \in \mathfrak{J}$ such that $\mathscr{D}^* = \mathscr{D}[\phi(\mathscr{D}^{**})/\mathscr{D}^{**}]$. The binary relation of $\mathscr{D}$ \emph{reducing} to $\mathscr{D}^*$ \emph{via} $\mathfrak{J}$---written $\mathscr{D} \leq_\mathfrak{J} \mathscr{D}^*$---is the reflexive-transitive closure of $\mathscr{D}$ immediately reducing to $\mathscr{D}^*$ relative to $\mathfrak{J}$.
\end{definition}

\begin{definition}
    An \emph{argument} is a pair $\langle \mathscr{D}, \mathfrak{J} \rangle$. It is \emph{closed} \emph{(}resp. \emph{open}\emph{)} iff $\mathscr{D}$ is so.
\end{definition}

\begin{definition}
    $\langle \mathscr{D}, \mathfrak{J} \rangle$ is \emph{valid} on $\mathfrak{B} \Longleftrightarrow$
    \begin{itemize}
        \item $\mathscr{D}$ is closed $\Longrightarrow$
        \begin{itemize}
            \item the conclusion of $\mathscr{D}$ is atomic $\Longrightarrow \mathscr{D} \leq_\mathfrak{J} \mathscr{D}^*$ with $\mathscr{D}^* = \langle T, f \rangle$ and $T \in \texttt{DER}_\mathfrak{B}$
            \item the conclusion of $\mathscr{D}$ is not atomic $\Longrightarrow \mathscr{D} \leq_\mathfrak{J} \mathscr{D}^*$ with $\mathscr{D}^*$ canonical whose immediate sub-structures are valid on $\mathfrak{B}$ when paired with $\mathfrak{J}$
        \end{itemize}
        \item $\mathscr{D}$ is open with undischarged assumptions $A_1, ..., A_n \Longrightarrow \forall \sigma, \forall A_i \ (i \leq n)$, if $\langle \sigma(A_i), \mathfrak{J}^+ \rangle$ is closed valid on $\mathfrak{B} \ (i \leq n)$, then $\langle \mathscr{D}^\sigma, \mathfrak{J}^+ \rangle$ is valid on $\mathfrak{B}$.
    \end{itemize}
\end{definition}
\noindent Argumental validity on an atomic base as defined here is non-monotonic, since there are $\langle \mathscr{D}, \mathfrak{J} \rangle$ and $\mathfrak{B}$ such that $\langle \mathscr{D}, \mathfrak{J} \rangle$ is valid on $\mathfrak{B}$ but, for some extension $\mathfrak{C}$ of $\mathfrak{B}$, $\langle \mathscr{D}, \mathfrak{J} \rangle$ is not valid on $\mathfrak{C}$. E.g., take $\mathfrak{J} = \mathfrak{B} = \emptyset, \mathfrak{C} = \{p\}$, and let $\mathscr{D}$ be the one-step argument structure from assumption $p$ to conclusion $q$.

It is easy to see that the open one-step structure consisting of just an application of $\vee_\lambda$ from $A \vee B$ to $A$ is valid with respect to $\{\phi_{\vee_\lambda}\}$ on $\mathfrak{B}$, as soon as there is no closed valid argument for $B$ on $\mathfrak{B}$. Suppose we are given a closed $\langle \mathscr{D}, \mathfrak{J} \rangle$ for $A \vee B$ valid on $\mathfrak{B}$. Then $\mathscr{D}$ reduces modulo $\mathfrak{J}$ to a closed canonical structure with conclusion $A \vee B$, whose immediate sub-structure $\mathscr{D}^*$ is closed valid relative to $\mathfrak{J}$ and $\mathfrak{B}$. The conclusion of $\mathscr{D}^*$ must then be $A$. So, when $\vee_\lambda$ is appended to $\mathscr{D}$, the closure $\mathscr{D}^\sigma$ thereby obtained of the one-step structure we started with, reduces modulo $\mathfrak{J}$ to the following argument structure.

\begin{prooftree}
    \AxiomC{$\mathscr{D}^*$}
    \noLine
    \UnaryInfC{$A$}
    \UnaryInfC{$A \vee B$}
    \RightLabel{$\vee_\lambda$}
    \UnaryInfC{$A$}
\end{prooftree}
\noindent In turn, the latter reduces to $\mathscr{D}^*$ by one application of $\phi_{\vee_\lambda}$. Since $\langle \mathscr{D}^*, \mathfrak{J} \rangle$ is valid on $\mathfrak{B}$, $\langle \mathscr{D}^*, \mathfrak{J} \cup \{\phi_{\vee_\lambda}\} \rangle$ is valid on $\mathfrak{B}$, hence $\langle \mathscr{D}^\sigma, \mathfrak{J} \cup \{\phi_{\vee_\lambda}\} \rangle$ is also so. The result follows by arbitrariness of $\langle \mathscr{D}, \mathfrak{J} \rangle$.

\begin{definition}
    $\langle \mathscr{D}, \mathfrak{J} \rangle$ is \emph{valid} iff $\langle \mathscr{D}, \mathfrak{J} \rangle$ is valid on every $\mathfrak{B}$.
\end{definition}

\noindent We can now prove that the open one-step structure consisting of just an application of $\texttt{Wk}$ from $A \rightarrow B$ to $(A \wedge C) \rightarrow B$ is valid with respect to $\{\phi_\wedge, \phi_{\texttt{Wk}}\}$. Take any $\mathfrak{B}$ and any closed $\langle \mathscr{D}, \mathfrak{J} \rangle$ for $A \rightarrow B$ valid on $\mathfrak{B}$. Then, $\mathscr{D}$ reduces modulo $\mathfrak{J}$ to a closed canonical structure with conclusion $A \rightarrow B$, whose immediate sub-structure $\mathscr{D}^*$ is open from $A$ to $B$ valid relative to $\mathfrak{J}$ and $\mathfrak{B}$. Therefore, when $\texttt{Wk}$ is appended to $\mathscr{D}$, the closure thereby obtained of the one-step structure we started with, reduces modulo $\mathfrak{J}$ to the following argument structure.

\begin{prooftree}
    \AxiomC{$[A]_1$}
    \noLine
    \UnaryInfC{$\mathscr{D}^*$}
    \noLine
    \UnaryInfC{$B$}
    \RightLabel{$1$}
    \UnaryInfC{$A \rightarrow B$}
    \LeftLabel{$\mathscr{D}^{**} = \ $}
    \RightLabel{$\texttt{Wk}$}
    \UnaryInfC{$(A \wedge C) \rightarrow B$}
\end{prooftree}
In turn, the latter reduces to the following by one application of $\phi_{\texttt{Wk}}$.

\begin{prooftree}
    \AxiomC{$[A \wedge C]_1$}
    \RightLabel{$\wedge_{E, 1}$}
    \UnaryInfC{$A$}
    \noLine
    \UnaryInfC{$\mathscr{D}^*$}
    \noLine
    \UnaryInfC{$B$}
    \LeftLabel{$\mathscr{D}^{***} = \ $}
    \RightLabel{$1$}
    \UnaryInfC{$(A \wedge C) \rightarrow B$}
\end{prooftree}
The open argument structure obtained by putting the given application of $\wedge_{E, i}$ on top of $\mathscr{D}^*$ is valid from $A \wedge C$ to $B$ relative to $\mathfrak{J} \cup \{\phi_\wedge\}$ and $\mathfrak{B}$. For, any closed structure for $A \wedge C$ valid relative to $(\mathfrak{J} \cup \{\phi_\wedge\})^+$ and $\mathfrak{B}$ reduces modulo $(\mathfrak{J} \cup \{\phi_\wedge\})^+$ to a canonical structure whose immediate sub-structures are valid relative to $(\mathfrak{J} \cup \{\phi_\wedge\})^+$ and $\mathfrak{B}$. Hence so is the left-hand immediate sub-structure $\mathscr{D}^{****}$ for $A$. The closure of the given open structure from $A \wedge C$ to $B$ obtained by replacing the open assumption $A \wedge C$ by the given closed structure for $A \wedge C$ will hence reduce modulo $(\mathfrak{J} \cup \{\phi_\wedge\})^+$---by applying $\phi_\wedge$---to a closure $(\mathscr{D}^*)^\sigma$ of $\mathscr{D}^*$ where the open assumption $A$ of $\mathscr{D}^*$ is replaced by $\mathscr{D}^{****}$. Since $\langle \mathscr{D}^*, \mathfrak{J} \rangle$ was valid on $\mathfrak{B}$, $\langle (\mathscr{D}^*)^\sigma, (\mathfrak{J} \cup \{\phi_\wedge\})^+ \rangle$ is valid on $\mathfrak{B}$. Given that the immediate sub-structure of $\mathscr{D}^{***}$ is valid relative to $\mathfrak{J} \cup \{\phi_\wedge\}$ on $\mathfrak{B}$, $\langle \mathscr{D}^{***}, \mathfrak{J} \cup \{\phi_\wedge\} \rangle$ is also so, hence $\langle \mathscr{D}^{***}, \mathfrak{J} \cup \{\phi_\wedge, \phi_{\texttt{Wk}}\} \rangle$ is valid on $\mathfrak{B}$, and so is hence $\langle \mathscr{D}^{**}, \mathfrak{J} \cup \{\phi_\wedge, \phi_{\texttt{Wk}}\} \rangle$ too. The result now follows by arbitrariness of $\langle \mathscr{D}, \mathfrak{J} \rangle$ and of $\mathfrak{B}$.

\begin{definition}
    $A$ is a \emph{consequence} of $\Gamma$ \emph{on} $\mathfrak{B}$ \emph{in the sense of reducibility semantics}---written $\Gamma \models^\alpha_\mathfrak{B} A$---iff there is $\langle \mathscr{D}, \mathfrak{J} \rangle$ from $\Gamma$ to $A$ which is valid on  $\mathfrak{B}$. $A$ is a \emph{logical consequence} of $\Gamma$ \emph{in the sense of reducibility semantics}---written $\Gamma \models^\alpha A$---iff there is a valid $\langle \mathscr{D}, \mathfrak{J} \rangle$ from $\Gamma$ to $A$.
\end{definition}
\noindent From what shown above, we hence have $A \vee B \models^\alpha_\mathfrak{B} A$ as soon as $\not\models^\alpha_\mathfrak{B} B$, and $A \rightarrow B \models^\alpha (A \wedge C) \rightarrow B$.

Observe that it is not \emph{a priori} clear whether logical consequence is tantamount to consequence on every atomic base. However, the following is trivially provable.

\begin{proposition}
    $\Gamma \models^\alpha A \Longrightarrow \forall \mathfrak{B} \ (\Gamma \models^\alpha_\mathfrak{B} A)$.
\end{proposition}
\noindent Whether we also have the inverse of this implication is in fact tantamount to the question whether logical consequence in reducibility semantics is fully equivalent to logical consequence in standard base semantics. I deal with this in the next section, where it is also shown that the issue is strictly interrelated with the question whether the inversion of the following implication holds.

\begin{proposition}
    $\Gamma \models^\alpha_\mathfrak{B} A \Longrightarrow \ (\models^\alpha_\mathfrak{B} \Gamma \Longrightarrow \ \models^\alpha_\mathfrak{B} A)$.
\end{proposition}

\begin{proof}
    Suppose $\Gamma \models^\alpha_\mathfrak{B} A$, i.e., that there is a
    \begin{prooftree}
        \AxiomC{$\Gamma$}
        \noLine
        \UnaryInfC{$\mathscr{D}$}
        \noLine
        \UnaryInfC{$A$}
    \end{prooftree}
    which is valid on $\mathfrak{B}$ when paired with $\mathfrak{J}$. Suppose $\models^\alpha_\mathfrak{B} \Gamma$ with $\Gamma = \{A_1, ..., A_n\}$, and take any closed $\langle \mathscr{D}_i, \mathfrak{J}^+ \rangle$ for $A_i \ (i \leq n)$ which is valid on $\mathfrak{B}$. By Definition 15,
    \begin{prooftree}
        \AxiomC{$\mathscr{D}_1, ..., \mathscr{D}_n$}
        \noLine
        \UnaryInfC{$\Gamma$}
        \noLine
        \UnaryInfC{$\mathscr{D}$}
        \noLine
        \UnaryInfC{$A$}
    \end{prooftree}
    is closed valid on $\mathfrak{B}$ when paired with $\mathfrak{J}^+$, whence $\models^\alpha_\mathfrak{B} A$.
\end{proof}

\section{Reducibility semantics \emph{vs} base semantics}

Let us now turn to the comparison between the three kinds of PTS introduced in Section 3. I start from showing that, under suitable conditions (whose limits will be discussed in the concluding remarks), reducibility semantics is equivalent to standard base semantics. The following preliminary results must be priorly proved, showing that Prawitz's reducibility semantics yields in fact a standard consequence relation over a ‘‘model", if (sets of) reductions are allowed to be ‘‘liberal" enough.

\begin{proposition}
The following facts hold:
    \begin{enumerate}
        \item $A \in \texttt{ATOM}_\mathscr{L} \Longrightarrow \ (\models^\alpha_\mathfrak{B} A \Longleftrightarrow \ \vdash_\mathfrak{B} A)$
        \item $A = B \wedge C \Longrightarrow \ (\models^\alpha_\mathfrak{B} A \Longleftrightarrow \ \models^\alpha_\mathfrak{B} B$ and $\models^\alpha_\mathfrak{B} C)$
        \item $A = B \vee C \Longrightarrow \ (\models^\alpha_\mathfrak{B} A \Longleftrightarrow \ \models^\alpha_\mathfrak{B} B$ or $\models^\alpha_\mathfrak{B} C)$
        \item $A = B \rightarrow C \Longrightarrow \ (\models^\alpha_\mathfrak{B} A \Longleftrightarrow B \models^\alpha_\mathfrak{B} C)$
        \item $(\models^\alpha_\mathfrak{B} \Gamma \Longrightarrow \ \models^\alpha_\mathfrak{B} A) \Longrightarrow \Gamma \models^\alpha_\mathfrak{B} A$
        \end{enumerate}
\end{proposition}

\begin{proof}
    The only non-trivial part is point 5. When $\Gamma = \{A_1, ..., A_n\}$, let $\tau = \langle \mathscr{D}_1, ..., \mathscr{D}_n \rangle$ be a sequence of closed argument structures for $A_1, ..., A_n$ respectively such that, for some $\mathfrak{J}_i$, $\langle \mathscr{D}_i, \mathfrak{J}_i \rangle$ is valid on $\mathfrak{B} \ (i \leq n)$. By assumption, there is then a closed argument structure for $A$, say $\mathscr{D}^\tau_A$, such that, for some $\mathfrak{J}^\tau_A$, $\langle \mathscr{D}^\tau_A, \mathfrak{J}^\tau_A \rangle$ is valid on $\mathfrak{B}$. Consider now the rule $\{\langle \langle \mathscr{D}_1, ..., \mathscr{D}_n \rangle, A \rangle \}$, to whose only element the closed argument structure

    \begin{prooftree}
        \AxiomC{$\mathscr{D}_1$}
        \noLine
        \UnaryInfC{$A_1$}
        \AxiomC{$\dots$}
        \AxiomC{$\mathscr{D}_n$}
        \noLine
        \UnaryInfC{$A_n$}
        \LeftLabel{$\mathscr{D} = \ $}
        \TrinaryInfC{$A$}
    \end{prooftree}
    is associated, and consider the mapping $\phi^\tau(\mathscr{D}) = \mathscr{D}^\tau_A$. Next, let $\mathbb{F}$ be the set of all the mappings which satisfy Definition 12. Then, we have $\phi^\tau \in \mathbb{F}$ since (a) both $\mathscr{D}$ and $\mathscr{D}^\tau_A$ are closed, so condition 1 in Definition 12 is satisfied, and (b) for every $\sigma$, $\mathscr{D}^\sigma = \mathscr{D}$ and $(\mathscr{D}^\tau_A)^\sigma = \mathscr{D}^\tau_A$, so $\phi^\tau(\mathscr{D}^\sigma) = \phi^\tau(\mathscr{D}) = \mathscr{D}^\tau_A = (\mathscr{D}^\tau_A)^\sigma = \phi^\tau(\mathscr{D})^\sigma$, so condition 2 in Definition 12 is satisfied too. But then, the open argument structure
    \begin{prooftree}
        \AxiomC{$A_1$}
        \AxiomC{$\dots$}
        \AxiomC{$A_n$}
        \LeftLabel{$\mathscr{D}^* =$}
        \TrinaryInfC{$A$}
    \end{prooftree}
    is valid on $\mathfrak{B}$ when $\mathbb{F}$ is its reduction-set. To see this, take any sequence $\tau = \langle \mathscr{D}_1, ..., \mathscr{D}_n \rangle$ of closed argument structures for $A_1, ..., A_n$ respectively such that, for some $\mathfrak{J}_i$, $\langle \mathscr{D}_i, \mathfrak{J}_i \rangle$ is valid on $\mathfrak{B}$, and let $\sigma(A_i) = \mathscr{D}_i \ (i \leq n)$. We must now prove that $\langle (\mathscr{D}^*)^\sigma, \mathbb{F} \cup \bigcup_{i \leq n} \mathfrak{J}_i \rangle$ is closed valid for $A$ on $\mathfrak{B}$ which, because $\mathfrak{J}_i \subseteq \mathbb{F} \ (i \leq n)$, can be simplified to showing that $\langle (\mathscr{D}^*)^\sigma, \mathbb{F} \rangle$ is closed valid for $A$ on $\mathfrak{B}$, i.e., $(\mathscr{D}^*)^\sigma$ reduces via $\mathbb{F}$ to a closed canonical argument structure for $A$ whose immediate sub-structures are valid on $\mathfrak{B}$ when paired with $\mathbb{F}$. But reducibility is transitive---see \cite{prawitz1973}---so we can limit ourselves to showing that $(\mathscr{D}^*)^\sigma$ reduces through $\mathbb{F}$ to a closed argument structure for $A$ which is valid on $\mathfrak{B}$ relative to $\mathbb{F}$. Now, since $\mathbb{F}$ contains a reduction $\phi^\tau$ as above, $(\mathscr{D}^*)^\sigma$ reduces via $\mathbb{F}$ to a closed $\mathscr{D}^\tau_A$ as above, namely, such that, for some $\mathfrak{J}^\tau_A$, $\langle \mathscr{D}^\tau_A, \mathfrak{J}^\tau_A \rangle$ is closed valid for $A$ on $\mathfrak{B}$---in particular, $\phi((\mathscr{D}^*)^\sigma) = (\mathscr{D}^*)^\sigma[\mathscr{D}^\tau_A/(\mathscr{D}^*)^\sigma]$. Given that $\mathfrak{J}^\tau_A \subseteq \mathbb{F}$, we have that $\langle \mathscr{D}^\tau_A, \mathbb{F} \rangle$ is also closed valid for $A$ on $\mathfrak{B}$.
\end{proof}
\noindent Some remarks are now in order concerning the proof of point 5 of Proposition 3. First of all, observe that the mapping $\phi^\tau$ used in the proof is just a ‘‘pointer" that sends a given sequence of closed argument structures for $A_1, ..., A_n$ onto some other closed argument structure for $A$ whose existence is guaranteed by the assumption $\models^\alpha_\mathfrak{B} \Gamma \Longrightarrow \ \models^\alpha_\mathfrak{B} A$. This ‘‘pointer" may of course not really ‘‘build" the closed argument structure for $A$ out of the given closed argument structures for $A_1, ..., A_n$. In fact, it might well be an instance on specific inputs of a constant function $\kappa$, say $\Gamma = \{p\}$ and $A = B \rightarrow B$, where we require that, for every closed argument structure $\mathscr{D}$ for $q$, given the argument structures

\begin{prooftree}
    \AxiomC{$\mathscr{D}$}
    \noLine
    \UnaryInfC{$q$}
    \LeftLabel{$\mathscr{D}^* = \ $}
    \UnaryInfC{$B \rightarrow B$}
    \AxiomC{$[B]$}
    \LeftLabel{$\mathscr{D}^{**} = \ $}
    \UnaryInfC{$B \rightarrow B$}
    \noLine
    \BinaryInfC{}
\end{prooftree}
we have $\kappa(\mathscr{D}^*) = \mathscr{D}^{**}$. In general, however, if we do not allow for classical logic in the meta-language, different sequences $\tau_1, \tau_2$ of closed argument structures for $A_1, ..., A_n$ may ‘‘point" to different closed argument structures for $A$, so that $\phi^{\tau_1} \neq \phi^{\tau_2}$.\footnote{Thus nothing changes if, instead of defining a $\phi^\tau$ for each argument structure obtained by appending the argument structure $\mathscr{D}^*$ in the proof to a sequence $\tau = \langle \mathscr{D}_1, ..., \mathscr{D}_n \rangle$ of closed argument structures for $A_1, ..., A_n$ respectively, we define a single $\phi$ for the rule consisting of all these argument structures in such a way that, if the instance $\iota$ of the rule is obtained by appending $\mathscr{D}^*$ to the sequence $\tau$, then $\phi(\iota) = \phi^\tau(\iota)$.} What is important for the theorem to go through is that $\mathbb{F}$ contains \emph{all} these ‘‘pointers". Clearly, one could choose the smaller set containing just the ‘‘pointers" plus the reduction-sets associated to the closed argument structures for $A_1, ..., A_n, A$. This would certainly be a more constructive choice, but still not an entirely constructive one. I say more about this in the Concluding Remarks. When on the contrary we allow for classical logic in the meta-language, then the situation is smoother, and in fact we need not pass through specific sequences of closed argument structures for $A_1, ..., A_n$---see \cite{piccolomininote}. For, if there is $A_i \ (i \leq n)$ such that there is no closed argument for $A_i$ valid on $\mathfrak{B}$, then the argument structure $\mathscr{D}^*$ in the proof above can be justified by the empty function, so by the empty reduction-set. Otherwise, there is some closed valid argument $\langle \mathscr{D}, \mathfrak{J} \rangle$ for $A$, so we can consider the rule $R$ consisting of all the argument structures obtained by appending $\mathscr{D}^*$ to arbitrary argument structures for $A_1, ..., A_n$, and define the following constant function $\kappa$, working as a reduction for $R$: for every instance $\iota$ of $R$, $\kappa(\iota) = \mathscr{D}$. Then, $\langle \mathscr{D}^*, \kappa \cup \mathfrak{J} \rangle$ is valid on $\mathfrak{B}$. I shall come back to this too in the Concluding Remarks.

Something more should be said about the constructive reading of the meta-logical constants as used in the proof of point 5 of Proposition 3---as remarked above, when the meta-logic is classical, there are no problems. The proof is based on the idea that the implication $\models^\alpha_\mathfrak{B} \Gamma \Longrightarrow \ \models^\alpha_\mathfrak{B} A$---where for simplicity we can assume $\Gamma = \{B\}$---yields a constructive function $\phi_1$ from proofs for

\begin{center}
    there is $\langle \mathscr{D}, \mathfrak{J} \rangle$ for $B$ valid on $\mathfrak{B}$
\end{center}
to proofs for

\begin{center}
    there is $\langle \mathscr{D}^*, \mathfrak{J}^* \rangle$ for $A$ valid on $\mathfrak{B}$
\end{center}
namely, from pairs $\langle \langle \mathscr{D}, \mathfrak{J} \rangle, \pi \rangle$, with $\pi$ (meta-)proof of the fact that $\langle \mathscr{D}, \mathfrak{J} \rangle$ is valid on $\mathfrak{B}$, to pairs $\langle \langle \mathscr{D}^*, \mathfrak{J}^* \rangle, \pi^* \rangle$, with $\pi^*$ (meta-)proof of the fact that $\langle \mathscr{D}^*, \mathfrak{J}^* \rangle$ is valid on $\mathfrak{B}$. The intuition is that one can extract from this a constructive function $\phi_2$ from proofs for $B$ to proofs for $A$, that is, from left-hand elements of pairs in the domain of $\phi_1$ to left-hand elements of pairs in the co-domain of $\phi_1$, and possibly---since reductions must go \emph{from argument structures to argument structures}---a further constructive function $\phi_3$ from left-hand elements of pairs in the domain of $\phi_2$ to left-hand elements of pairs in the co-domain of $\phi_2$. Now, if one assumes (as I have implicitly done) that, for each $\langle \mathscr{D}, \mathfrak{J} \rangle$ for $B$, the fact that $\langle \mathscr{D}, \mathfrak{J} \rangle$ is valid on $\mathfrak{B}$ has \emph{at most} one (meta-)proof $\texttt{P}(\langle \mathscr{D}, \mathfrak{J} \rangle)$, $\phi_2$ is readily found, i.e.,

\begin{center}
    $\phi_2(\langle \mathscr{D}, \mathfrak{J} \rangle) \stackrel{\text{def}}{=} \texttt{lp}(\phi_1(\langle \langle \mathscr{D}, \mathfrak{J} \rangle, \texttt{P}(\langle \mathscr{D}, \mathfrak{J} \rangle) \rangle))$
\end{center}
where $\texttt{lp}$ is the standard left-projection on a pair. On the other hand, it seems unnatural to assume that each $\mathscr{D}$ for $B$ is valid on $\mathfrak{B}$ with respect to \emph{at most} one $\mathfrak{J}$. So, it is unclear how to attain a definition of $\phi_3$. Be it as it may, the set of all the mappings which satisfy Definition 12, i.e., $\mathbb{F}$ above, must contain a mapping from any $\mathscr{D}$ for $B$ such that $\langle \mathscr{D}, \mathfrak{J} \rangle$ is valid on $\mathfrak{B}$ for some $\mathfrak{J}$, to the image of $\langle \mathscr{D}, \mathfrak{J} \rangle$ under $\phi_2$.

Definition 4 and Proposition 3 yield the equivalence between standard base semantics and reducibility semantics relative to consequence over a base.

\begin{theorem}
    $\Gamma \models_\mathfrak{B} A \Longleftrightarrow \Gamma \models^\alpha_\mathfrak{B} A$.
\end{theorem}

\begin{proof}
    By induction when $\Gamma = \emptyset$, and by reduction to the case when $\Gamma \neq \emptyset$ via point 5 of Proposition 3 and the corresponding clause of Definition 4.
\end{proof}

\noindent The equivalence between the two approaches at issue relative to \emph{logical} consequence is instead obtained by combining Theorem 1 with Proposition 1, plus the following. 

\begin{proposition}
    $\forall \mathfrak{B} \ (\Gamma \models^\alpha_\mathfrak{B} A) \Longrightarrow \Gamma \models^\alpha A$.
\end{proposition}

\begin{proof}

The proof is similar to that of point 5 of Proposition 3. Assume $\Gamma = \{A_1, ..., A_n\}$. We know that, for every $\mathfrak{B}$, we have an argument structure

\begin{prooftree}
    \AxiomC{$A_1$}
    \AxiomC{$\dots$}
    \AxiomC{$A_n$}
    \noLine
    \TrinaryInfC{$\mathscr{D}_\mathfrak{B}$}
    \noLine
    \UnaryInfC{$A$}
\end{prooftree}
such that $\langle \mathscr{D}_\mathfrak{B}, \mathfrak{J}_\mathfrak{B} \rangle$ is valid on $\mathfrak{B}$ for some $\mathfrak{J}_\mathfrak{B}$. Hence, for all $\mathfrak{B}$, the set $\mathbb{F}$ of all the mappings satisfying Definition 12 contains $\phi_\mathfrak{B}$ which, given closed $\mathscr{D}^i_\mathfrak{B}$-s with conclusion $A_i$ such that, for some $\mathfrak{J}^i_\mathfrak{B}$, $\langle \mathscr{D}^i_\mathfrak{B}, \mathfrak{J}^i_\mathfrak{B} \rangle$ is valid on $\mathfrak{B} \ (i \leq n)$, maps the argument structure

\begin{prooftree}
    \AxiomC{$\mathscr{D}^1_\mathfrak{B}$}
    \noLine
    \UnaryInfC{$A_1$}
    \AxiomC{$\dots$}
    \AxiomC{$\mathscr{D}^n_\mathfrak{B}$}
    \noLine
    \UnaryInfC{$A_n$}
    \LeftLabel{$\mathscr{D}^\sigma = \ \ \ $}
    \TrinaryInfC{$A$}
\end{prooftree}
onto the argument structure
\begin{prooftree}
    \AxiomC{$\mathscr{D}^1_\mathfrak{B}$}
    \noLine
    \UnaryInfC{$A_1$}
    \AxiomC{$\dots$}
    \AxiomC{$\mathscr{D}^n_\mathfrak{B}$}
    \noLine
    \UnaryInfC{$A_n$}
    \noLine
    \LeftLabel{$\mathscr{D}^*_\mathfrak{B} = \ \ \ $}
    \TrinaryInfC{$\mathscr{D}_\mathfrak{B}$}
    \noLine
    \UnaryInfC{$A$}
\end{prooftree}
We want to prove that $\langle \mathscr{D}^\sigma, \mathbb{F} \cup \bigcup_{i \leq n} \mathfrak{J}^i_\mathfrak{B} \rangle$ is valid on $\mathfrak{B}$ which, since $\mathbb{F} \cup \bigcup_{i \leq n} \mathfrak{J}^i_\mathfrak{B} = \mathbb{F}$, can be reduced to showing that $\langle \mathscr{D}^\sigma, \mathbb{F} \rangle$ is valid on $\mathfrak{B}$. $\mathscr{D}^\sigma$ reduces via $\phi_\mathfrak{B} \in \mathbb{F}$ to $\mathscr{D}^*_\mathfrak{B}$ which, by assumption, is such that $\langle \mathscr{D}^*_\mathfrak{B}, \mathfrak{J}_\mathfrak{B} \cup \bigcup_{i \leq n} \mathfrak{J}^i_\mathfrak{B} \rangle$ is valid on $\mathfrak{B}$. But $\mathfrak{J}_\mathfrak{B} \cup \bigcup_{i \leq n} \mathfrak{J}^i_\mathfrak{B} \subseteq \mathbb{F}$, thus $\langle \mathscr{D}^*_\mathfrak{B}, \mathbb{F} \rangle$ is also valid on $\mathfrak{B}$.
\end{proof}

\noindent All the remarks about the proof of point 5 in Proposition 3 apply to the proof of Proposition 4 too---the constructive reading of the meta-logical constants requires now that we take atomic bases as parameters, see \cite{piccolomininote}.\footnote{The advantage of adopting classical logic in the meta-language  in the case of proof of point 5 in Proposition 3---i.e., the possibility of defining a single constant function rather than a class of input-depending ‘‘pointers"---is of no use here since, even with classical logic in the meta-language, a constant function is now required for each atomic base, and there is \emph{a priori} no guarantee that these constant functions are identical to one another.} Then, we easily obtain the following.

\begin{theorem}
    $\Gamma \models A \Longleftrightarrow \Gamma \models^\alpha A$.
\end{theorem}

\begin{proof}
    By Definition 6, $\Gamma \models A \Longleftrightarrow \forall \mathfrak{B} \ (\Gamma \models_\mathfrak{B} A)$. By Theorem 1, $\forall \mathfrak{B} \ (\Gamma \models_\mathfrak{B} A) \Longleftrightarrow \forall \mathfrak{B} \ (\Gamma \models^\alpha_\mathfrak{B} A)$. By Propositions 1 and 4, $\forall \mathfrak{B} \ (\Gamma \models^\alpha_\mathfrak{B} A) \Longleftrightarrow \Gamma \models^\alpha A$.
\end{proof}

\noindent For what said at the end of Section 3.1, i.e., that $\texttt{IL}$ is incomplete over $\models$ (with classical meta-logic), from Theorem 2 we then obtain the following.

\begin{theorem}
    $\texttt{IL}$ is incomplete over $\models^\alpha$ (with classical meta-logic).
\end{theorem}

Before moving to the next topic, let me comment a little bit more upon the results I have been focusing on in this section. However easy their proofs are (and under the constraint that consequence is always taken to hold between a formula and a \emph{finite} set of formulas), Theorem 1 and Theorem 2 provide a seemingly interesting link between standard base semantics and reducibility semantics. As already hinted at in the Introduction, the latter can be seen as a sort of ‘‘realizer-interpretation" of the former---see \cite{GheorghiuPym}---in the sense that, whenever $\Gamma \models_\mathfrak{B} A$ or $\Gamma \models A$ hold, a ‘‘witness" can be found in terms of an open argument from $\Gamma$ to $A$ which is valid on $\mathfrak{B}$, i.e., $\Gamma \models^\alpha_\mathfrak{B} A$, or logically valid, i.e., $\Gamma \models^\alpha A$, respectively. Conversely, standard base semantics might be understood as a sort of type-abstraction from reducibility semantics, meaning that, whenever we are given an open argument ‘‘witnessing" either $\Gamma \models^\alpha_\mathfrak{B} A$ or $\Gamma \models^\alpha A$, we disregard \emph{how} the connection between $\Gamma$ and $A$ is established, and just look at the input and output types of it, i.e., $\Gamma \models_\mathfrak{B} A$ and $\Gamma \models A$ respectively.

Theorem 1 and Theorem 2 confirm this rough intuition---thus implying incompleteness of $\texttt{IL}$ over $\models^\alpha$ as per Theorem 3. If one thinks about that, however, this is far from trivial. To my mind, the distinction between reducibility semantics and standard base semantics is best understood if, \emph{mutatis mutandis}, we look at how truth is dealt with in an another major constructivist framework, i.e., Martin-L\"{o}f's intuitionistic type theory \cite{martin-loef}. Leaving atomic bases aside, let us concede that closed valid arguments for $A$, or open valid arguments for $A$ from $\Gamma = \{A_1, ..., A_n\}$, can be read in terms of proof-objects involved in type-theoretic judgements of the form

\begin{center}
    $a : A$ \\
    
    $b(x_1, ..., x_n) : A \ (x_1 : A_1, ..., x_n : A_n)$
\end{center}
respectively. The first judgement says that $a$ is a proof-object of type $A$---namely, of proposition $A$, given the Curry-Howard isomorphism and the underlying \emph{formulas-as-types conception} \cite{howard}. The second judgement says instead that $b(a_1, ..., a_n)$ is a proof-object of type $A$, whenever $a_i$ is a proof-object of type $A_i \ (i \leq n)$. In intuitionistic type theory, we also have judgements of the form

\begin{center}
    $A \ \texttt{true}$ \\

    $A \ \texttt{true} \ (A_1 \ \texttt{true}, ..., A_n \ \texttt{true})$
\end{center}
where the second judgement says that $A$ is true whenever $A_i$ is true ($i \leq n$), and where $\texttt{true}$ behaves according to the following rules respectively---in the right-hand rule I allow myself for a type-theoretic notational abuse, which is however harmless in the present context:
\begin{prooftree}
    \AxiomC{$a : A$}
    \RightLabel{$\texttt{T}$}
    \UnaryInfC{$A \ \texttt{true}$}
    \AxiomC{$b(x_1, ..., x_n) : A \ (x_1 : A, ..., x_n : A_n)$}
    \RightLabel{$\texttt{DT}$}
    \UnaryInfC{$A \ \texttt{true} \ (A_1 \ \texttt{true}, ..., A_n \ \texttt{true})$}
    \noLine
    \BinaryInfC{}
\end{prooftree}
Thus, in Martin-L\"{o}f's intuitionistic type theory, $\texttt{true}$ is always explained in terms of a given proof-object, which can be exhibited as a truth-maker of the proposition whose truth is being asserted. This is similar to what happens in reducibility semantics, where the consequence relation is not primitive, but defined as existence of a suitable valid argument. In standard base semantics, we start instead from the consequence relation outright, without requiring something more primitive. Sticking to the comparison with Martin-L\"{o}f's intuitionistic type theory, it is thus as if we started from a notion of truth which is not explained in terms of the rules $\texttt{T}$ and $\texttt{DT}$ above, say
\begin{center}
    $A \ \texttt{true}^*$ \\

    $A \ \texttt{true}^* \ (A_1 \ \texttt{true}^*, ..., A_n \ \texttt{true}^*)$.
\end{center}
If one accepts this reconstruction, the question about the relation between reducibility semantics and standard base semantics can be thus understood in rough type-theoretic terms as the question whether at least one of the following rules is acceptable:
\begin{prooftree}
    \AxiomC{$A \ \texttt{true} \ (A_1 \ \texttt{true}, ..., A_n \ \texttt{true})$}
    \UnaryInfC{$A \ \texttt{true}^* \ (A_1 \ \texttt{true}^*, ..., A_n \ \texttt{true}^*)$}
    \AxiomC{$A \ \texttt{true}^* \ (A_1 \ \texttt{true}^*, ..., A_n \ \texttt{true}^*)$}
    \UnaryInfC{$A \ \texttt{true} \ (A_1 \ \texttt{true}, ..., A_n \ \texttt{true})$}
    \noLine
    \BinaryInfC{}
\end{prooftree}
hence whether at least one of the following rules is acceptable:
\begin{prooftree}
    \AxiomC{$b(x_1, ..., x_n) : A \ (x_1 : A_1, ..., x_n : A_n)$}
    \UnaryInfC{$A \ \texttt{true}^* \ (A_1 \ \texttt{true}^*, ..., A_n \ \texttt{true}^*)$}
    \AxiomC{$A \ \texttt{true}^* \ (A_1 \ \texttt{true}^*, ..., A_n \ \texttt{true}^*)$}
    \UnaryInfC{$b(x_1, ..., x_n) : A \ (x_1 : A_1, ..., x_n : A_n)$}
    \noLine
    \BinaryInfC{}
\end{prooftree}
with $n \geq 0$. As far as this goes, Theorem 1 and Theorem 2 say that, to a large extent, the questions can be answered positively.

Let us now turn to a quick comparison between $\models^\alpha$ and $\models^s$---remarking that, since $\models^\alpha$ has been proved to be equivalent to $\models$, the results I shall present here also hold when we replace $\models^\alpha$ with $\models$.

Because of the fact that Sanqvist's base semantics explains disjunction in an elimination like way, we will not have full equivalence here, but a weaker result. For seeing this, given an arbitrary $\mathfrak{B}$, consider the statement

\begin{equation}
    \forall \Gamma \ \forall A \ (\Gamma \models^s_\mathfrak{B} A \Longrightarrow \Gamma \models^\alpha_\mathfrak{B} A).
\end{equation}
Then, we can prove the following.

\begin{theorem}
    $\forall \mathfrak{B} \ ((1) \Longrightarrow \forall \Gamma \ \forall A \ (\Gamma \models^\alpha_\mathfrak{B} A \Longrightarrow \Gamma \models^s_\mathfrak{B} A))$.
\end{theorem}

\begin{proof}
    By induction on $A$ when $\Gamma = \emptyset$, and by reducing to the closed case---via point 5 of Proposition 3---when $\Gamma \neq \emptyset$. 
\end{proof}

Thus, contrarily to what happens with the comparison between $\models$ and $\models^\alpha$, a ‘‘local" comparison is not possible between $\models^\alpha$ and $\models^s$---by ‘‘local" I mean that, for \emph{each} $\Gamma$ and $A$, $\Gamma \models^\alpha_\mathfrak{B} A$ implies $\Gamma \models^s_\alpha A$, and vice versa. The only comparability we can hope for is a ‘‘conditionally global" one, i.e., \emph{if} we can go \emph{everywhere} from Sandqvist to Prawitz, \emph{then} we can also go \emph{everywhere} from Prawitz to Sandqvist.

Of course, the ‘‘conditionally global" comparability of $\models^\alpha$ and $\models^s$ may not be without value. For example, it may shed light on mutual properties of $\models^\alpha$ and $\models^s$ \emph{at the level of logical validity}, or it may provide us with insights concerning the relative soundness and completeness of given logics over one or the other approach. This is what happens in the monotonic framework, after introducing notions of ‘‘point-wise" soundness and completeness. But now, significant differences between monotonic and non-monotonic PTS emerge, as I point out in the Section 5.

\paragraph{Remark on the connection with (in)completeness results} In the next section, I will deal with the difference between monotonic and non-monotonic PTS relative to the comparability results of this section---as well as relative to some notions of ‘‘point-wise" soundness and completeness which I will introduce in due course. But the distinction between monotonic and non-monotonic PTS relative to the results discussed in this section can be also discussed with reference to the issue of the completeness or incompleteness of given logics over standard or Sandqvist's base semantics, and reducibility semantics. To my knowledge, the only concrete result of this kind that we have in the non-monotonic framework is Piecha and Schroeder-Heister's proof that $\texttt{IL}$ is incomplete over $\models$ \cite{piechaschroeder-heister2019}, from which we drew Theorem 3, so $\texttt{IL}$ is incomplete over $\models^\alpha$ too---relevant results for a (partly elimination-based) approach are to be found also in \cite{takemura}. It is instead unsettled which logic is \emph{complete} over $\models$ (hence, as per Theorem 2, over $\models^\alpha$), and this also applies to Sandqvist's non-monotonic base semantics, although in this case a proof of incompleteness of $\texttt{IL}$ is missing.\footnote{One should specify (as I will do in the Concluding Remarks too) that we are here taking the PTS-notions of consequence as defined relative to atomic bases \emph{with no upper bound} on the complexity level. But these notions can be defined also relative to atomic bases whose complexity level \emph{does have} an upper bound, say, consequence over bases whose level is \emph{at most} some positive number $n$. In Piecha and Schroeder-Heister's incompleteness proof it is essential to require that the complexity be unlimited, since the proof is based on a function which translates sets $\Gamma$ of disjunction-free formulas into atomic bases $\Gamma^\circ$: left-nested implications in elements of $\Gamma$ will normally make the level of $\Gamma^\circ$ increase, and since there might be as many left-nested implication as one may want, no complexity constraint can be put on atomic bases. But of course, one might well have (in)completeness results which depend on the uppermost complexity level of the atomic bases, say a given logic is complete over one of the three PTS notions of consequence at issue in this paper, \emph{provided the bases} have level at most (or at least) some fixed valued $n$, and so on. While we do have something similar in monotonic PTS (see below), we again have nothing similar in the non-monotonic context.}

In monotonic PTS, the situation is much richer instead. As proved in \cite{piechaschroeder-heisterdecampossanz, piechaschroeder-heister2019}, $\texttt{IL}$ is generally incomplete over monotonic standard base semantics (and over monotonic reducibility semantics, since a result similar to Theorem 2 can be proved in monotonic PTS too).\footnote{The expression ‘‘generally" means here that the result holds both when the complexity of atomic bases is allowed to be unlimited, and when it is instead limited by some upper bound.} Also, Stafford proved in \cite{stafford1} that inquisitive logic is complete over $\models$ when the complexity of the atomic bases is $\geq 2$ and, under the same circumstances, Sandqvist proved as said in $\cite{sandqvist}$ that $\texttt{IL}$ is complete over his own variant of base semantics.

As shown in \cite{piccolominiinversion}, the completeness and incompleteness results for monotonic PTS can be combined with the monotonic variant of Theorem 2 (and with notions of ‘‘point-wise" soundness and completeness mentioned below) to obtain non-trivial insights relative to such issues as whether the three kinds of monotonic PTS at issue in this paper are or not base-comparable---where base-comparability roughly means that, when $A$ is a consequence on $\Gamma$ on a given base over one of these three approaches, then $A$ is also a consequence of $\Gamma$ on that base in another approach. But due to the lack of suitable completeness or incompleteness results, no such use of Theorem 4 seems to be available in the case of non-monotonic PTS (the same holds for the notions of ‘‘point-wise" soundness and completeness, but for different reasons). This does not mean, however, that future completeness or incompleteness results for non-monotonic PTS might show consequences of Theorem 4 other than that of providing a connection between $\models^s$ and $\models^\alpha$ (and, via Theorem 2, with $\models^s$ and $\models$).

\section{On base-soundness and base-completeness}

Let me now turn to a comparison between the comparability results established in Section 4, and analogous results which can be proved for the monotonic variants of $\models, \models^\alpha$ and $\models^s$. In particular, I aim at showing in this section how some notions of ‘‘point-wise" soundness and completeness behave in non-monotonic PTS, and at drawing from this certain conclusions on the differences between the monotonic and the non-monotonic approaches.

The general point is the following: ‘‘point-wise" soundness and completeness have similar effects in both the frameworks---in particular, $\texttt{IL}$ is not base-complete over any of the kinds of PTS at issue in this paper---though to a \emph{limited extent} only. More specifically, while these notions can be used in the monotonic picture to obtain some positive information about the connection between standard base semantics and reducibility semantics, on the one hand, and Sandqvist's base semantics on the other, they will have no such effect in the non-monotonic context. And the reason of this difference is in turn that, relative to the ‘‘point-wise" notions, some general results can be proved in non-monotonic PTS which, due to the monotonicity constraint, will fail in the monotonic variant.

Let us write $\models_M, \models^\alpha_M$ and $\models^s_M$ the \emph{monotonic} variants of $\models, \models^\alpha$ and $\models^s$ respectively---these obtain from Definitions 4, 5 and 15, and from Proposition 3, by requiring consequence over $\mathfrak{B}$ when $\Gamma \neq \emptyset$ to range over arbitrary extensions of $\mathfrak{B}$. Now, a result in all ways identical to Theorem 4 can be proved for these monotonic variants too---where one has of course to bring extensions of atomic bases in, namely, the global invertibility from Sandqvist to Prawitz and vice versa holds for any extension of a given atomic base. In turn, this can be used in combination with notions of base-soundness and base-completeness of given logics over PTS, to single out a sufficient condition for $\models^\alpha_M$ (or $\models_M$) and $\models^s_M$ to be equivalent at the level of logical validity, as well as for a logic to be complete over $\models^\alpha_M$ (or $\models_M$).

Base-completeness is defined as follows---in the definition, $\Vdash$ stands for $\models^\alpha_{(M)}, \models_{(M)}$, or $\models^s_{(M)}$, while $\Sigma$ is a super-intuitionistic logic, i.e., a recursive set of super-intuitionistic rules with a derivability relation $\vdash_\Sigma$, possibly implementing an atomic base $\mathfrak{B}$, in which case derivability is written $\vdash_{\Sigma \cup \mathfrak{B}}$---see \cite{piccolominiinversion} for precise definitions.

\begin{definition}
    $\Sigma$ is \emph{base-complete} over $\Vdash$ iff $\forall \Gamma \ \forall A \ \forall \mathfrak{B} \ (\Gamma \Vdash_\mathfrak{B} A \Longrightarrow \Gamma \vdash_{\Sigma \cup \mathfrak{B}} A)$.
\end{definition}

\noindent Completeness is clearly the same as base-completeness, but without reference to atomic bases.

\begin{definition}
    $\Sigma$ is \emph{complete} over $\Vdash$ iff $\forall \Gamma \ \forall A \ (\Gamma \Vdash A \Longrightarrow \Gamma \vdash_\Sigma A)$.
\end{definition}
\noindent From this one easily proves what follows.

\begin{proposition}
    $\Sigma$ base-complete over $\Vdash \ \Longrightarrow \ \Sigma$ complete over $\Vdash$.
\end{proposition}

\begin{proof}
    Assume $\Gamma \Vdash A$. So $\forall \mathfrak{B} \ (\Gamma \Vdash_\mathfrak{B} A)$. Instantiate this on the empty base $\mathfrak{B}^\emptyset$, i.e., $\Gamma \Vdash_{\mathfrak{B}^\emptyset} A$. By base-completeness, $\Gamma \vdash_{\Sigma \cup \mathfrak{B}^\emptyset} A$ but, since $\mathfrak{B}^\emptyset$ is empty, we have $\Sigma \cup \mathfrak{B}^\emptyset = \Sigma$, so $\Gamma \vdash_\Sigma A$.
\end{proof}

\noindent As for base-soundness, we have instead what follows---the proof of Proposition 6 in the monotonic context would be slightly different, see \cite{piccolominiinversion}.

\begin{definition}
    $\Sigma$ is \emph{base-sound} over $\Vdash$ iff $\forall \Gamma \ \forall A \ \forall \mathfrak{B} \ (\Gamma \vdash_{\Sigma \cup \mathfrak{B}} A \Longrightarrow \Gamma \Vdash_\mathfrak{B} A)$.
\end{definition}

\begin{definition}
    $\Sigma$ is \emph{sound} over $\Vdash$ iff $\forall \Gamma \ \forall A \ (\Gamma \vdash_\Sigma A \Longrightarrow \Gamma \Vdash A)$.
\end{definition}

\noindent If we assume that $\Sigma$ is monotonic, i.e., for every extension $\Omega$ of $\Sigma$, $\Gamma \vdash_\Sigma A \Longleftrightarrow \Gamma \vdash_\Omega A$, we can prove the following.

\begin{proposition}
    $\Sigma$ base-sound over $\Vdash \ \Longrightarrow \Sigma$ sound over $\Vdash$.
\end{proposition}

\begin{proof}
    If $\Gamma \vdash_\Sigma A$, by monotonicity, $\forall \mathfrak{B} \ (\Gamma \vdash_{\Sigma \cup \mathfrak{B}} A)$. By base-soundness, $\forall \mathfrak{B} \ (\Gamma \Vdash_\mathfrak{B} A)$, i.e., $\Gamma \Vdash A$.
\end{proof}

\begin{proposition}
    $\texttt{IL}$ is base-sound over $\Vdash$.
\end{proposition}

\begin{corollary}
    $\texttt{IL}$ is sound over $\Vdash$ (see e.g. \cite{prawitz1973}).
\end{corollary}

Now, the monotonic versions of Theorem 4, Definitions 18, 19, 20 and 21, Proposition 5, 6 and 7, and Corollary 1, jointly yield that $\models^\alpha_M$ (or $\models_M$) and $\models^s_M$ are equivalent at the level of logical validity if there is a $\Sigma$ which is base-complete over $\models^s_M$ and base-sound over $\models^\alpha_M$ (or $\models_M$). The same condition is sufficient for completeness over $\models^\alpha_M$ (or $\models_M$)---see \cite{piccolominiinversion}. These facts, though, cannot yield any ‘‘positive" information about $\texttt{IL}$ over $\models^\alpha_M$ (or $\models_M$), as $\texttt{IL}$ can be shown \emph{not to be} base-complete over monotonic $\Vdash$ for any (either limited or unlimited) upper bound on the complexity of atomic bases. This is proved by taking into account what de Campos Sanz, Piecha and Schroeder-Heister called the \emph{export principle} \cite{piechaschroeder-heisterdecampossanz}, the latter being a general method for shifting from consequence over an atomic base $\mathfrak{B}$ to logical consequence and back, via a translation of $\mathfrak{B}$ into a set of disjunction-free formulas $\mathfrak{B}^*$, thus $\Gamma \Vdash_{\mathfrak{B}} A \Longleftrightarrow \Gamma, \mathfrak{B}^* \Vdash A$.

The idea is now asking whether a given logic $\Sigma$ on atomic base $\mathfrak{B}$ enjoys the export principle, provided one is allowed to assume or discharge (through rules in $\mathfrak{B}$), not only sets of formulas $\Gamma$, but also sets of atomic rules $\mathfrak{C}$, i.e., $\Gamma, \mathfrak{C} \vdash_{\Sigma \cup \mathfrak{B}} A \Longleftrightarrow \Gamma, \mathfrak{C}^*, \mathfrak{B}^* \vdash_\Sigma A$. $\texttt{IL}$ enjoys this property, and this implies that base-completeness of $\texttt{IL}$ over monotonic $\Vdash$ (with an either limited or unlimited upper bound on the complexity of atomic bases) is tantamount to completeness over $\Vdash$ plus $\Vdash$ enjoying the export principle (relative to the given bound on the complexity of atomic bases). Which means that base-completeness of $\texttt{IL}$ over monotonic $\Vdash$ is always a void notion, as Piecha and Schroeder-Heister proved in general that export principle plus completeness of $\texttt{IL}$ over $\Vdash$ imply incompleteness of $\texttt{IL}$ over $\Vdash$ \cite{piechaschroeder-heister2019}. This extends to any $\Sigma$ which enjoys the export principle and does not derive Harrop rule---the latter being Piecha and Schroeder-Heister's counterexample to completeness of $\texttt{IL}$ over $\models_M$.

The same results would hold, \emph{mutatis mutandis}, in the non-monotonic context too. The situation would be here even smoother, because one can immediately observe that non-monotonic $\Vdash$ \emph{does not} enjoy the export principle.\footnote{The translation $^*$ from atomic bases to sets of disjunction free formulas and back is defined by induction on the level of the atomic rules:

\begin{itemize}
    \item $R$ is of level $0 \Longrightarrow R^* = R$
    \item $R$ is of level $k$, hence has the form
    \begin{prooftree}
        \AxiomC{$[\mathfrak{C}_1]$}
        \noLine
        \UnaryInfC{$A_1$}
        \AxiomC{$\dots$}
        \AxiomC{$[\mathfrak{C}_n]$}
        \noLine
        \UnaryInfC{$A_n$}
        \TrinaryInfC{$A$}
    \end{prooftree}
    where each
    \begin{prooftree}
        \AxiomC{$\mathfrak{C}_i$}
        \noLine
        \UnaryInfC{$A_i$}
    \end{prooftree}
    is a rule $R_i$ of level at most $k - 1 \Longrightarrow R^* = \bigwedge_{i \leq n} R^*_i \rightarrow A$.
\end{itemize}
Then, when $\mathfrak{B} = \{R_1, ..., R_n\}$, we put $\mathfrak{B}^* = \{R_1^*, ..., R_n^*\}$. Based on the above given definition of what it means for $\Vdash$ to enjoy the export principle, we see immediately what follows.

\begin{fact}

Non-monotonic $\Vdash$ does not enjoy the export principle
\end{fact}
\noindent For otherwise, from $p \Vdash_{\mathfrak{B}^\emptyset} q$, we would have $p, (\mathfrak{B}^\emptyset)^* \Vdash q$, i.e., $p \Vdash q$, which is clearly false.} Thus, this would extend from $\texttt{IL}$ to any logic which enjoys the export principle, whether it derives Harrop rule or not---as the full use of Piecha and Schroeder-Heister's incompleteness theorem is no longer needed.

However, one needs not go through all this, for now, a major difference between the monotonic and the non-monotonic approaches occurs. This depends on the following result---below, we say that $\Sigma$ is \emph{closed under substitutions} iff for every $\Gamma, A, A_1, ..., A_n$, if the propositional variables of $\Gamma \cup \{A\}$ are $p_1, ..., p_n$, then

\begin{center} 
$\Gamma \vdash_\Sigma A \Longleftrightarrow \Gamma[A_1, ..., A_n/p_1, ..., p_n] \vdash_\Sigma A[A_1, ..., A_n/p_1, ..., p_n]$
\end{center}
and $\Vdash$ is again restricted to $\models, \models^\alpha$ and $\models^s$.

\begin{proposition}
    $\Sigma$ is consistent and closed under substitutions $\Longrightarrow \Sigma$ not base-complete over $\Vdash$.
\end{proposition}

\begin{proof}
    Suppose that $\Sigma$ is consistent, closed under substitutions, and base-complete over $\Vdash$. By base-completeness, $p \Vdash_{\mathfrak{B}^\emptyset} q \Longrightarrow p \vdash_{\Sigma \cup \mathfrak{B}^\emptyset} q$. But $\Sigma \cup \mathfrak{B}^\emptyset = \Sigma$, so that $p \Vdash_{\mathfrak{B}^\emptyset} q \Longrightarrow p \vdash_\Sigma q$. As the antecedent of this implication holds, we infer $p \vdash_\Sigma q$. By closure under substitutions, hence, for every $A$ and $B$, $A \vdash_\Sigma B$, which cannot be by the consistency assumption. Hence, $\Sigma$ is not base-complete over $\Vdash$.
\end{proof}

\begin{corollary}
    $\texttt{IL}$ is not base-complete over $\Vdash$.\footnote{Let me remark that the base-incompleteness of $\texttt{IL}$ over $\models$ and $\models^\alpha$ follows easily, with classical meta-logic, from Theorems 2 and 3, plus Proposition 5. Anyway, Corollary 2 has an independent interest since, first, it holds constructively and, second, it holds also for Sandqvist's variant.}
\end{corollary}

\noindent Proposition 8 and Corollary 2 jointly show that the conditions for base-incompleteness of arbitrary logics over non-monotonic proof-theoretic consequence are much weaker than those required for the monotonic case, and that base-incompleteness of $\texttt{IL}$ in the non-monotonic framework obtains in a much smoother way than over the monotonic counter-part. But besides that, they also seem to suggest that the result about the sufficient condition for equivalence between $\models^\alpha_M$ (or $\models_M$) and $\models^s_M$ and for completeness over $\models^\alpha_M$ (or $\models_M$)---which I mentioned above, i.e., existence of a $\Sigma$ which is base-complete over $\models^s_M$ and base-sound over $\models^\alpha_M$ (or $\models_M$)---is quite trivial in the non-monotonic case. For, if $\Sigma$ is consistent and closed under substitutions, the implication is vacuously true as $\Sigma$ cannot be base-complete over $\models^s$, and if $\Sigma$ is inconsistent, the implication is again vacuously true, since $\Sigma$ cannot in this case be base-sound over $\models^\alpha$ (or $\models$), because $\not\models^\alpha \bot$ (resp. $\not\models \bot$). Thus, one might think that the only case in which such a condition could be informative is the one where $\Sigma$ is not closed under substitutions. But this is false, due to the following result---where $\Vdash^1$ and $\Vdash^2$ are either $\models$, or $\models^\alpha$, or $\models^s$, and are not necessarily distinct.

\begin{proposition}
    $\Sigma$ sound over $\Vdash^1 \ \Longrightarrow \Sigma$ not base-complete over $\Vdash^2$.
\end{proposition}

\begin{proof}
    Suppose $\Sigma$ base-complete over $\Vdash^2$. Thus, for all $\mathfrak{B}$, $p \Vdash^2_\mathfrak{B} q \Longrightarrow p \vdash_{\Sigma \cup \mathfrak{B}} q$. Since $p \Vdash^2_{\mathfrak{B}^\emptyset} q$, we have $p \vdash_{\Sigma \cup \mathfrak{B}^\emptyset} q$, i.e., $p \vdash_\Sigma q$. Since $\Sigma$ is sound over $\Vdash^1$, we have $p \Vdash^1 q$, which is false---consider the base which consists of axiom $p$ only. Thus, $\Sigma$ cannot be base-complete over $\Vdash^2$.
\end{proof}

\begin{corollary}
    $\Sigma$ base-sound over $\Vdash^1 \ \Longrightarrow \Sigma$ not base-complete over $\Vdash^2$.
\end{corollary}

\noindent Observe that base-incompleteness of $\texttt{IL}$ over non-monotonic $\Vdash$ also stems from Proposition 9 plus Corollary 1---respectively, Corollary 3 plus Proposition 7. To conclude, I note that both Proposition 8 and Corollary 2, and Proposition 9 and Corollary 3, as well as the consequences they imply on the comparison between monotonic and non-monotonic PTS, obviously depend on the non-monotonicity of the consequence relations at issue. For, in the \emph{monotonic} context, we do not have $p \Vdash_{\mathfrak{B}^\emptyset} q$, as now the consequence relation must hold for arbitrary extensions of $\mathfrak{B}^\emptyset$, including that consisting of just axiom $p$.

\section{Concluding remarks}

To conclude, I just want to observe first of all that the results I presented so far also hold if we assume that atomic bases have \emph{limited} complexity---say, when the consequence relation is restricted to atomic bases whose complexity is \emph{at most} $n$. Albeit I put no restriction of this kind, and allowed atomic bases to contain rules of \emph{any} complexity, this is immaterial as concerns the proposed comparison between $\models, \models^\alpha$ and $\models^s$.

Next, let me deal quickly with some limitations which the results presented in this paper may undergo when further constraints are put on the definitions of reducibility semantics.

The reduction set $\mathbb{F}$ used for justifying the equivalence between reducibility semantics and standard base semantics---i.e., Theorem 1 and Theorem 2 respectively---may be looked upon as not constructive enough. More specifically, $\mathbb{F}$ is obviously not liable to a description as, say, a recursive set of recursive rewriting rules for argument structures. One may appeal to a smaller set, say $\texttt{S}$, of the functions $\phi_\mathfrak{B}$ such that, given any closed $\mathscr{D}_i$ for $A_i$ such that, for some $\mathfrak{J}_i$, $\langle \mathscr{D}_i, \mathfrak{J}_i \rangle$ is valid on $\mathfrak{B}$, map instances $\mathscr{D}^\sigma$ with $\sigma(A_i) = \mathscr{D}_i \ (i \leq n)$ of

\begin{prooftree}
    \AxiomC{$A_1$}
    \AxiomC{}
    \noLine
    \UnaryInfC{$\dots$}
    \AxiomC{$A_n$}
    \LeftLabel{$\mathscr{D} =$}
    \TrinaryInfC{$A$}
\end{prooftree}
onto a closed $\mathscr{D}_A$ for $A$ such that $\langle \mathscr{D}_A, \mathfrak{J}_A \rangle$ is also valid on $\mathfrak{B}$ for some $\mathfrak{J}_A$. However, although any such $\phi_\mathfrak{B}$ may be looked upon as constructive, $\texttt{S}$ may again not be expressible as a recursive set of recursive rewriting rules. This seems to hold even when, using classical logic in the meta-language, one observes that $\phi_\mathfrak{B}$ is either a ‘‘pointer" which maps every $\mathscr{D}^\sigma$ onto a fixed closed $\mathscr{D}_A$ for $A$---if every $A_i$ is provable on $\mathfrak{B} \ (i \leq n)$---or the empty function $\phi_\emptyset$---if there is $A_i$ not provable on $\mathfrak{B} \ (i \leq n)$---i.e., for every $\sigma$,
\begin{center}
    $\phi_\mathfrak{B}(\mathscr{D}^\sigma) = \begin{cases}
        \mathscr{D}_A & \text{if} \ \forall i \leq n \ (\models^\alpha_\mathfrak{B} A_i) \\ \phi_\emptyset & \text{otherwise}
    \end{cases}$
\end{center}
---for more on this, see \cite{piccolomininote} while, for a general discussion about what a good reduction is, see \cite{ayhangoodreductions}. Thus, under a stricter reading of what reductions and sets of reductions are to be, which yields a stricter reading of reducibility semantics, indicated by $\models^\alpha_S$, it is not clear how to attain point 5 of Proposition 3 and Proposition 4, hence Theorem 1, Theorem 2 and Proposition 5 too---for a first discussion on this, see \cite{piccolominithesis, piccolominibook, piccolomininote}. In fact, more constructive examples of reductions for justifying non-intuitionistically derivable rules can be provided, again by adapting to a Prawitzian framework with argument structures and reductions the results proved by de Campos Sanz, Piecha and Schroeder-Heister (both in the monotonic and in the non-monotonic context)---see \cite{piechadecampossanz, piechaschroeder-heisterdecampossanz, piechaschroeder-heister2019}---or possibly via \emph{Pezlar's selector} for the Split rule---see \cite{pezlarselector}. These issues can be dealt with in future works---for a preliminary discussion, see again \cite{piccolomininote}.

Be that as it may, even with a stricter notion of reduction, it is still possible to prove an inversion result between $\models^\alpha_S$ and $\models$ similar to the one proved for $\models^\alpha$ and $\models^s$. Take any arbitrary $\mathfrak{B}$ and consider the statement
\begin{equation}
    \forall \Gamma \ \forall A \ (\Gamma \models_\mathfrak{B} A \Longrightarrow \Gamma \models^\alpha_{S, \mathfrak{B}} A).
\end{equation}
\noindent Then we have what follows.

\begin{theorem}
    $\forall \mathfrak{B} \ ((2) \Longrightarrow \forall \Gamma \ \forall A \ (\Gamma \models^\alpha_{S, \mathfrak{B}} A \Longrightarrow \Gamma \models_\mathfrak{B} A))$.
\end{theorem}

\noindent The proof is in all ways similar to that of Theorem 4.

\paragraph{Acknowledgments} I am very grateful to Ansten Klev, Hermógenes Oliveira, Thomas Piecha, Dag Prawitz, Antje Rumberg, Peter Schroeder-Heister, Will Stafford, and the anonymous reviewers, for precious observations and remarks which helped me improve previous versions of this paper. This work has been supported by the grant PI 1965/1-1 for the DFG project \emph{Revolutions and paradigms in logic. The case of proof-theoretic semantics}.

\bibliographystyle{abbrv}
\bibliography{bibliography}

\begin{thebibliography}{10}

\bibitem{ayhangoodreductions}
S.~Ayhan.
\newblock What are acceptable reductions? {Perspectives} from proof-theoretic semantics and type theory.
\newblock {\em Australasian journal of logic}, 20(3):412--428, 2023.
\newblock \url{https://doi.org/10.26686/ajl.v20i3.7692}.

\bibitem{piechadecampossanz}
W.~{de Campos Sanz} and T.~Piecha.
\newblock A critical remark on the {BHK} interpretation of implication.
\newblock {\em Philosophia Scientiae}, 18(3):13--22, 2014.
\newblock \url{https://doi.org/10.4000/philosophiascientiae.965}.

\bibitem{dummett1991}
M.~Dummett.
\newblock {\em The logical basis of metaphysics}.
\newblock Harvard University Press, 1993.

\bibitem{francez}
N.~Francez.
\newblock {\em Proof-theoretic semantics}.
\newblock College Publications, 2015.

\bibitem{gentzenuntersuchungen}
G.~Gentzen.
\newblock Untersuchungen \"{u}ber das logische {Schließen} i, ii.
\newblock {\em Mathematische Zeitschrift}, 39:176--210, 405--431, 1935.
\newblock \url{https://doi.org/10.1007/BF01201353}.

\bibitem{gentzen}
G.~Gentzen.
\newblock Untersuchungen über das logische {Schließen} i, ii.
\newblock {\em Mathematische Zeitschrift}, 39:176–210, 405–431, 1935.
\newblock \url{https://doi.org/10.1007/BF01201353}.

\bibitem{GheorghiuPym}
A.~Gheorghiu and D.~Pym.
\newblock From proof-theoretic validity to base-extension semantics for intuitionistic propositional logic.
\newblock \url{arXiv:2210.05344}, 2022.

\bibitem{howard}
W.~Howard.
\newblock The formulae-as-types notion of construction.
\newblock In J.~Seldin and R.~Hindley, editors, {\em To {H.B. Curry}: Essays on combinatory logic, lambda calculus and formalism}, pages 479--490. Academy Press, 1980.

\bibitem{martin-loef}
P.~{Martin-L\"{o}f}.
\newblock {\em Intuitionistic type theory}.
\newblock Bibliopolis, 1984.

\bibitem{hermogenespragmatist}
H.~Oliveira.
\newblock On {Dummett's} pragmatist justification procedure.
\newblock {\em Erkenntnis}, 86:429--455, 2021.
\newblock \url{https://doi.org/10.1007/s10670-019-00112-7}.

\bibitem{pezlarselector}
I.~Pezlar.
\newblock Constructive validity of a generalised {Kreisel-Putnam} rule.
\newblock {\em Studia Logica}, 2024.
\newblock \url{https://doi.org/10.1007/s11225-024-10129-x}.

\bibitem{piccolominithesis}
A.~{Piccolomini d'Aragona}.
\newblock {\em Dag Prawitz's theory of grounds}.
\newblock PhD Dissertation, Aix-Marseille University, Sapienza University of Rome, HAL Id: tel-02482320, 2019.

\bibitem{piccolominibook}
A.~{Piccolomini d'Aragona}.
\newblock {\em Prawitz's epistemic grounding. An investigation into the power of deduction}.
\newblock Synthese Library, Springer, 2023.

\bibitem{ptssquare}
A.~{Piccolomini d'Aragona}.
\newblock The proof-theoretic square.
\newblock {\em Synthese}, 201(219), 2023.
\newblock \url{https://doi.org/10.1007/s11229-023-04203-5}.

\bibitem{piccolomininote}
A.~{Piccolomini d'Aragona}.
\newblock A note on schematic validity and completeness in {Prawitz}'s semantics.
\newblock In F.~Bianchini, V.~Fano, and P.~Graziani, editors, {\em Current topics in logic and the philosophy of science. Papers from SILFS 2022 postgraduate conference}. College Publications, 2024.

\bibitem{piccolominiinversion}
A.~{Piccolomini d'Aragona}.
\newblock A comparison of three kinds of monotonic proof-theoretic semantics and the base-incompleteness of intuitionistic logic.
\newblock Submitted.

\bibitem{piechaschroeder-heisterdecampossanz}
T.~Piecha, W.~{de Campos Sanz}, and P.~Schroder-Heister.
\newblock Failure of completeness in proof-theoretic semantics.
\newblock {\em Journal of philosophical logic}, 44:321--335, 2015.
\newblock \url{https://doi.org/10.1007/s10992-014-9322-x}.

\bibitem{piechaschroeder-heisterbases}
T.~Piecha and P.~{Schroeder-Heister}.
\newblock Atomic systems in proof-theoretic semantics: two approaches.
\newblock In J.~Redmon, O.~{Pombo Martins}, and A.~{Nepomuceno Fernández}, editors, {\em Epistemology, knowledge and the impact of interaction}, pages 47--62. Springer, 2016.
\newblock \url{https://doi.org/10.1007/978-3-319-26506-3}.

\bibitem{piechaschroeder-heister2019}
T.~Piecha and P.~{Schroeder-Heister}.
\newblock Incompleteness of intuitionistic propositional logic with respect to proof-theoretic semantics.
\newblock {\em Studia Logica}, 107(1):47--62, 2019.
\newblock \url{https://doi.org/10.1007/s11225-018-9823-7}.

\bibitem{prawitz1965}
D.~Prawitz.
\newblock {\em Natural deduction. A proof-theoretical study}.
\newblock Almqvist \& Wiskell, 1965.

\bibitem{prawitz1971}
D.~Prawitz.
\newblock Ideas and results in proof theory.
\newblock In J.~E. Fenstad, editor, {\em Proceedings of the second Scandinavian logic symposium}, pages 235--307. Elsevier, 1971.
\newblock \url{https://doi.org/10.1016/S0049-237X(08)70849-8}.

\bibitem{prawitz1973}
D.~Prawitz.
\newblock Towards a foundation of a general proof-theory.
\newblock In P.~Suppes, L.~Henkin, A.~Joja, and G.~C. Moisil, editors, {\em Proceedings of the Fourth International Congress for Logic, Methodology and Philosophy of Science, Bucharest, 1971}, pages 225--250. Elsevier, 1973.
\newblock \url{https://doi.org/10.1016/S0049-237X(09)70361-1}.

\bibitem{sandqvist}
T.~Sandqvist.
\newblock Base-extension semantics for intuitionistic sentential logic.
\newblock {\em Logic journal of the IGPL}, 23(5):719--731, 2015.
\newblock \url{https://doi.org/10.1093/jigpal/jzv021}.

\bibitem{schroeder-heister2006}
P.~Schroeder-Heister.
\newblock Validity concepts in proof-theoretic semantics.
\newblock {\em Synthese}, 148:525--571, 2006.
\newblock \url{https://doi.org/10.1007/s11229-004-6296-1}.

\bibitem{schroeder-heisterSE}
P.~Schroeder-Heister.
\newblock Proof-theoretic semantics.
\newblock In E.~N. Zalta, editor, {\em The Stanford Encyclopedia of Philosophy}. 2018.

\bibitem{schroederheisterrolf}
P.~Schroeder-Heister.
\newblock Prawitz's completeness conjecture: a reassessment.
\newblock {\em Theoria}, 90(5):515--527, 2024.
\newblock \url{https://doi.org/10.1111/theo.12541}.

\bibitem{stafford1}
W.~Stafford.
\newblock Proof-theoretic semantics and inquisitive logic.
\newblock {\em Journal of philosophical logic}, 50:1199--1229, 2021.
\newblock \url{https://doi.org/10.1007/s10992-021-09596-7}.

\bibitem{stafford2}
W.~Stafford and V.~Nascimento.
\newblock Following all the rules: intuitionistic completeness for generalised proof-theoretic validity.
\newblock {\em Analysis}, 2023.
\newblock \url{https://doi.org/10.1093/analys/anac100}.

\bibitem{takemura}
R.~Takemura.
\newblock Investigation of {Prawitz's} completeness conjecture in phase semantic framework.
\newblock {\em Journal of Humanities and Sciences Nihon University}, 23(1):1--19, 2017.

\bibitem{wittgenstein}
L.~Wittgenstein.
\newblock {\em Philosophical investigations}.
\newblock Blackwell, 1953.

\end{thebibliography}

\end{document}